\documentclass[english,12pt,reqno,twoside]{amsart}

\usepackage{amssymb,amsmath,amsthm,soul,color}
\usepackage{graphicx}
\usepackage[noadjust]{cite}
\usepackage{t1enc}
\usepackage[utf8]{inputenc}
\usepackage{braket}
\usepackage[T1]{fontenc}
\usepackage{amsfonts}
\usepackage{a4,indentfirst,latexsym,graphics}
\usepackage[english]{babel}
\usepackage{mathtools}
\usepackage{fixmath}

\usepackage[left=2.1cm,right=2.1cm,top=2cm,bottom=2cm]{geometry}

\usepackage{url}
\usepackage{mathrsfs}

\usepackage{enumitem}
\usepackage[colorlinks=true,urlcolor=blue,
citecolor=red,linkcolor=blue,linktocpage,pdfpagelabels,
bookmarksnumbered,bookmarksopen]{hyperref}
\usepackage[hyperpageref]{backref}
\usepackage[colorinlistoftodos]{todonotes}

\usepackage[normalem]{ulem}
\usepackage{lmodern}

\makeatletter
\renewcommand{\l@subsection}{%
	\@tocline{2}{0pt}{\dimexpr1pc+2em\relax}{5pc}{}}
\renewcommand{\l@subsubsection}{%
	\@tocline{3}{0pt}{\dimexpr1pc+4em\relax}{7pc}{}}
\makeatother

\newtheorem{thm}{Theorem}[section]
\newtheorem{cor}[thm]{Corollary}
\newtheorem{lem}[thm]{Lemma}

\theoremstyle{definition}

\theoremstyle{remark}
\newtheorem{rem}[thm]{Remark}

\numberwithin{equation}{section}

\DeclareMathSymbol{\C}{\mathalpha}{AMSb}{"43}

\newcommand{\N}{\mathbb N}
\DeclareMathOperator{\Tr}{Tr}

\begin{document}
\title[Finite-rank fractional Lieb--Thirring inequality]{Optimizers for the finite-rank fractional Lieb--Thirring inequality}

\author[L. Chen]{Lu Chen}
\address{L. Chen 
	\newline \indent Key Laboratory of Algebraic Lie Theory and Analysis of Ministry of Education, School of Mathematics and Statistics, 
	\newline \indent Beijing Institute of Technology, 
	\newline \indent Haidian, 100081 Beijing, People's Republic of China}
\email{chenlu5818804@163.com}

\author[L. Yang]{Lianfeng Yang}
\address{L. Yang
	\newline \indent School of Mathematics and Statistics, 
	\newline \indent Beijing Institute of Technology, 
	\newline \indent Haidian, 100081 Beijing, People's Republic of China}
\email{yanglianfeng2021@163.com}

\date{}

\begin{abstract}
The finite-rank fractional Lieb--Thirring inequality provides an
estimate on the Riesz sum of the first $N$ negative eigenvalues of
the fractional Schr\"odinger operator $(-\Delta)^s-V$ in terms of
the $L^{\gamma+\frac{d}{2s}}(\mathbb R^d)$ norm of the potential $V$, where the Riesz exponent $\gamma>\max\left\{0,1-d/(2s)\right\}$.
We prove the existence of optimizers for every finite rank by
establishing a result on the spectrum in the case of dichotomy.
We further characterize the spectral gap and qualitative properties of optimizers, including the Euler--Lagrange equation, regularity, and asymptotic behavior at infinity.
In particular, under the condition $\gamma+{d}/{(2s)}\geq 3+{4s}/{d}$, we show that all normalized maximizing sequences are precompact modulo translations.
Moreover, we derive a strict binding inequality for the optimal constants in the finite-rank setting, which implies that the optimal fractional Lieb--Thirring constant, if it exists, cannot be attained by a potential with finitely many negative eigenvalues.
\end{abstract}

\maketitle
\tableofcontents

\section{Introduction}
The Lieb--Thirring inequalities provide bounds on the Riesz sums of the negative eigenvalues of Schr\"odinger operators $-\Delta-V$ in terms of the $L^p$ norms of the potentials $V$. They play a fundamental role in proving the stability of nonrelativistic matter \cite{LT75,LT76}.
Motivated by this work and the problem of stability of matter
in a relativistic setting, Daubechies \cite{Dau83} studied
Schr\"odinger operators with ``generalized'' kinetic energies, including fractional Schr\"odinger operators $(-\Delta)^s+V$.
Unlike the classical Laplacian case, the standard proof of
the stability of matter is no longer applicable here.
By establishing bounds on their negative eigenvalues, she eventually obtained lower bounds on the kinetic energy of fermions,
providing essential estimates for the future study of the
relativistic stability problem.
More recently, Frank, Lieb, and Seiringer \cite{FLS08} investigated Hardy--Lieb--Thirring inequalities for fractional Schr\"odinger operators with magnetic fields. As an application, they provided a  proof of the stability of relativistic matter, allowing the inclusion of arbitrary magnetic fields up to the critical nuclear charge.

Over the past five decades, these inequalities have been
extensively studied, resulting in a rich theory and numerous
important developments. 
Nevertheless, many interesting and challenging mathematical questions remain open \cite{Fra21}, including those concerning Lieb--Thirring inequalities for powers of the Laplacian.
Among them, the exact values of the optimal constants and
the existence of optimizers are two central questions.
We refer interested readers to
\cite{Fra09,FHJN21,LW00,Sch22,Hun07,LS10} and the references therein
for further background and related results.

As mentioned above, motivated by the importance of these inequalities, we study the fractional Lieb--Thirring inequality
(see \cite{Dau83,Fra18,Fra21,Nam23})
\begin{equation}\label{eq:fractional-LT}
	\sum_{j\geq1}|\lambda_j((-\Delta)^s-V)|^\gamma
\leq L_{\gamma,d,s}
\int_{\mathbb R^d}V_+^{\gamma+\frac{d}{2s}}\,dx
\end{equation}
for $0<s<1$, where
\[
\begin{cases}
	\gamma\geq 1-{d}/{(2s)}, & d<2s,\\[6pt]
	\gamma>0, & d=2s,\\[6pt]
	\gamma\geq0, & d>2s,
\end{cases}
\]
$V\in L^{\gamma+\frac{d}{2s}}(\mathbb R^d)$ is a real-valued function, $V_+=\max\{V,0\}$, and $\lambda_j((-\Delta)^s-V)$ denote the negative eigenvalues counted with multiplicity. 
A natural question is whether there exists an optimal potential and whether $L_{\gamma,d,s}$ can be attained by a potential with finitely many negative eigenvalues.

In the recent works \cite{FGL21,FGL21a,FGL25},
the authors proposed a new perspective on the optimal constants
in Lieb--Thirring inequalities by investigating the
finite-rank inequalities, where only the first $N$ eigenvalues
are considered. Inspired by these works, we provide an affirmative answer to the above question in the
finite-rank setting and show that any optimizer attaining
$L_{\gamma,d,s}$ in \eqref{eq:fractional-LT}, if it exists, must
have infinitely many negative eigenvalues. For related problems, see also \cite{CDGW26,CGW26,DLNN25,WZZ25}.

To this end, we consider the following
finite-rank fractional Lieb--Thirring inequality
\begin{equation}\label{eq:finite-rank-LT}
	\sum_{j=1}^N
	\left|\lambda_j((-\Delta)^s-V)\right|^\gamma
	\leq
	L^{(N)}_{\gamma,d,s}\int_{\mathbb R^d}V_+^{\gamma+\frac{d}{2s}}\,dx
\end{equation}
for $N\in\mathbb N$. Since $(-\Delta)^s-V\geq(-\Delta)^s-V_+$, we may assume without loss of generality that $V\geq0$ in the sequel.  We then define the optimal constant in \eqref{eq:finite-rank-LT} by
\begin{equation}\label{eq:constant}
	L^{(N)}_{\gamma,d,s}
	=
	\sup_{\substack{
			0\leq V\in L^{\gamma+\frac{d}{2s}}(\mathbb R^d)\\
			V\neq0}}
	\frac{\displaystyle
		\sum_{j=1}^N
		\left|\lambda_j((-\Delta)^s-V)\right|^\gamma}
	{\|V\|_{\gamma+\frac{d}{2s}}^{\gamma+\frac{d}{2s}}}.
\end{equation}
For a real-valued potential $V\in L^{\gamma+\frac{d}{2s}}(\mathbb R^d)$, the quadratic form
\[
q_V[u]
=
\|(-\Delta)^{s/2}u\|_2^2
-\int_{\mathbb R^d}V(x)|u(x)|^2\,dx,
\qquad u\in H^s(\mathbb R^d),
\]
defines the self-adjoint operator $H_V=(-\Delta)^s-V$
with $\sigma_{\mathrm{ess}}(H_V)=[0,\infty)$ (see Lemma~\ref{lem:form-compactness} for more details).
We denote by $\lambda_j(H_V)$ the $j$-th min--max level of
$H_V$, defined by
\[
\lambda_j(H_V)
:=
\inf_{\substack{
		\phi_1,\ldots,\phi_j\in H^s(\mathbb R^d)\\
		\langle\phi_k,\phi_\ell\rangle=\delta_{k\ell}
}}
\max_{\substack{
		u\in\operatorname{span}\{\phi_1,\ldots,\phi_j\}\\
		\|u\|_2=1
}}
q_V[u].
\]
By the min--max principle \cite{RSIV78}, it coincides with the $j$-th negative eigenvalue counted with multiplicity whenever it exists, and equals $0$ otherwise.

Observe that the variational problem in \eqref{eq:constant} is
invariant under the scaling
\[
V_t=t^{2s}V(t\cdot),\qquad t>0.
\]
Specifically, the map
$u\mapsto t^{d/2}u(t\cdot)$ is unitary on $L^2(\mathbb R^d)$
and maps $H^s(\mathbb R^d)$ onto itself.
Moreover, the scaling property of the fractional Laplacian
and a change of variables give
$q_{V_t}[t^{d/2}u(t\cdot)]=t^{2s}q_V[u]$.
By the min--max principle, we obtain
\begin{equation*}
	\lambda_j(H_{V_t})=t^{2s}\lambda_j(H_V),
	\qquad
	\|V_t\|_{\gamma+\frac{d}{2s}}^{\gamma+\frac{d}{2s}}
	=t^{2s\gamma}
	\|V\|_{\gamma+\frac{d}{2s}}^{\gamma+\frac{d}{2s}} .
\end{equation*}
Choosing $t=\|V\|_{\gamma+\frac{d}{2s}}^{-1/(2s\gamma)}$ yields
$\|V_t\|_{\gamma+\frac{d}{2s}}^{\gamma+\frac{d}{2s}}=1$ while preserving the quotient. Therefore,
\eqref{eq:constant} admits the equivalent form
\begin{equation}\label{eq:normalized-constant}
	L^{(N)}_{\gamma,d,s}
	=
	\sup_{\substack{
			0\leq V\in L^{\gamma+\frac{d}{2s}}(\mathbb R^d)\\
			\|V\|_{\gamma+\frac{d}{2s}}^{\gamma+\frac{d}{2s}}=1}}
	\sum_{j=1}^N|\lambda_j(H_V)|^\gamma .
\end{equation}
Throughout the paper, we focus on
\eqref{eq:normalized-constant} and refer to maximizers and
maximizing sequences for \eqref{eq:normalized-constant} as
optimizers and normalized maximizing sequences, respectively. Unless otherwise stated, we also assume
\[
d\in\N=\{1,2,\ldots\},\qquad 0<s<1,\qquad \gamma>\max\left\{0,1-d/(2s)\right\},\qquad p=\gamma+d/{(2s)}.
\]
Moreover, it is easy to check that $L^{(N)}_{\gamma,d,s}>0$ by testing with a smooth compactly supported potential. Together with the fractional Lieb--Thirring inequality \eqref{eq:fractional-LT}, this gives
\begin{equation}\label{eq:monotonicity}
	0<L^{(1)}_{\gamma,d,s}\leq L^{(2)}_{\gamma,d,s}\leq\cdots \leq L^{(N)}_{\gamma,d,s}
	\le L_{\gamma,d,s}<\infty.
\end{equation}

We are now in a position to state our main results. Our first result concerns the existence of optimizers for the finite-rank fractional Lieb--Thirring inequality.

\begin{thm}\label{thm:existence}
	Let $d\in\mathbb N$, $0<s<1$, and $\gamma>\max\left\{0,1-\frac{d}{2s}\right\}$.
	Then $L^{(N)}_{\gamma,d,s}$ admits an optimizer for every $N\in\mathbb N$.
\end{thm}

\begin{rem}\label{rem:binding-compactness}
	For $N\geq2$, assume in addition that
	\begin{equation}\label{eq:binding-criterion}
		\max_{1\leq k\leq\lfloor N/2\rfloor}
		L^{(k)}_{\gamma,d,s}
		<
		L^{(N)}_{\gamma,d,s}.
	\end{equation}
	Then every normalized maximizing sequence for
	$L^{(N)}_{\gamma,d,s}$ is precompact modulo translations in
	$L^p(\mathbb R^d)$. In particular, dichotomy in the proof of
	Theorem~\ref{thm:existence} is excluded.
	More precisely, recalling the proof of Theorem~\ref{thm:existence}, if dichotomy
	occurred, then \eqref{eq:plateau-split} would yield
	$L^{(k)}_{\gamma,d,s}
	=
	L^{(N-k)}_{\gamma,d,s}
	=
	L^{(N)}_{\gamma,d,s}$
	for some $1\leq k\leq N-1$. Noting that $1\leq\min\{k,N-k\}\leq\lfloor N/2\rfloor$, we have
	$L^{(\min\{k,N-k\})}_{\gamma,d,s}=L^{(N)}_{\gamma,d,s}$, and hence
	$$\max_{1\leq k\leq\lfloor N/2\rfloor}L^{(k)}_{\gamma,d,s}
	\geq L^{(\min\{k,N-k\})}_{\gamma,d,s}=L^{(N)}_{\gamma,d,s},$$
	which contradicts \eqref{eq:binding-criterion}.
\end{rem}

Our second result describes the qualitative properties of optimizers, including the spectral gap, the Euler--Lagrange equation, regularity, and asymptotic behavior at infinity.
We also note that the resolvent kernel estimates used below in \eqref{eq:kernel-coefficient} are known;
see \cite{FLS16}.

\begin{thm}\label{thm:structure}
Let $d,N\in\mathbb N$, $0<s<1$, and
$\gamma>\max\left\{0,1-{d}/{(2s)}\right\}$,
with $p=\gamma+{d}/{(2s)}$.
Assume that $V\geq0$ attains $L^{(N)}_{\gamma,d,s}$ and
$\|V\|_p^p=1$.
For brevity, write $\lambda_j:=\lambda_j(H_V)$ for $j\geq1$, and
denote by $K$ the number of strictly negative eigenvalues among
the first $N$ min--max levels, counted with multiplicity.
Then the following assertions hold.	
\begin{enumerate}[label=\textup{(\roman*)},ref=\textup{(\roman*)}]
\item\label{item:structure-EL}
We have $V>0$ almost everywhere in $\mathbb R^d$ and
\begin{equation}\label{eq:gap}
	\lambda_K<\lambda_{K+1}.
\end{equation}
Moreover, there exist real orthonormal eigenfunctions
$u_1,\ldots,u_K$ of $H_V$ corresponding to
$\lambda_1,\ldots,\lambda_K$, respectively,
such that the following identity holds almost everywhere
in $\mathbb R^d$:
\begin{equation}\label{eq:EL}
	V(x)=
	\left(\frac{\gamma}{pL^{(N)}_{\gamma,d,s}}
	\sum_{j=1}^K|\lambda_j|^{\gamma-1}|u_j(x)|^2
	\right)^{\frac{1}{p-1}}.
\end{equation}
These eigenfunctions also satisfy the following system
in the weak sense:
\begin{equation}\label{eq:EL-equation}
	(-\Delta)^su_j-
	\left(\frac{\gamma}{pL^{(N)}_{\gamma,d,s}}
	\sum_{\ell=1}^K|\lambda_\ell|^{\gamma-1}|u_\ell|^2
	\right)^{\frac{1}{p-1}}u_j
	=\lambda_j u_j,
	\qquad 1\leq j\leq K.
\end{equation}	

\item\label{item:structure-tail}
The eigenfunctions satisfy $u_j\in H^{2s}(\mathbb R^d)\cap C_0(\mathbb R^d)$ for $1\leq j\leq K$, and $V$ admits a continuous representative in
$C_0(\mathbb R^d)$.
Moreover, $u_1>0$ and $V>0$ everywhere.

For $a>0$, denote by $G_a$ the resolvent kernel
of $((-\Delta)^s+a)^{-1}$.
Let $a_{d,s}>0$ be the constant determined by
\begin{equation}\label{eq:kernel-coefficient}
	\lim_{|x|\to\infty}|x|^{d+2s}G_a(x)
	=a_{d,s}a^{-2},\qquad a>0.
\end{equation}
Then, as $|x|\to\infty$, we also have
\begin{equation}\label{eq:state-tail}
	u_j(x)=a_{d,s}|\lambda_j|^{-2}
	\left(\int_{\mathbb R^d}V(y)u_j(y)\,dy\right)
	|x|^{-(d+2s)}
	+o\left(|x|^{-(d+2s)}\right),
	\qquad 1\leq j\leq K,
\end{equation}
and
\begin{equation}\label{eq:potential-tail}
	V(x)
	=
	\left[
	\frac{\gamma}{pL^{(N)}_{\gamma,d,s}}
	\sum_{j=1}^K|\lambda_j|^{\gamma-1}
	\left(
	a_{d,s}|\lambda_j|^{-2}
	\int_{\mathbb R^d}V(y)u_j(y)\,dy
	\right)^2
	\right]^{\frac{1}{p-1}}
	|x|^{-\frac{2(d+2s)}{p-1}}+
	o\left(|x|^{-\frac{2(d+2s)}{p-1}}\right).
\end{equation}
\end{enumerate}		
\end{thm}

\begin{rem}\label{rem:L1-threshold}
Observe that $	a_{d,s}|\lambda_1|^{-2}
\int_{\mathbb R^d}V(y)u_1(y)\,dy>0$ since  $V>0$ and $u_1>0$ everywhere. So the leading coefficient in
\eqref{eq:potential-tail} is strictly positive. This yields the integrability criterion
\begin{equation}\label{eq:L1-threshold}
	V\in L^1(\mathbb R^d)
	\quad\Longleftrightarrow\quad
	p<3+\frac{4s}{d}.
\end{equation}
Indeed, since $V\in C_0(\mathbb R^d)$,
it is locally integrable.
On the other hand, by applying \eqref{eq:potential-tail}, we deduce that $V(x)\asymp |x|^{-\frac{2(d+2s)}{p-1}}$
as $|x|\to\infty$.
Consequently, for sufficiently large $R>0$,
\[
\int_{|x|>R}V(x)\,dx<\infty
\quad\Longleftrightarrow\quad
\int_R^\infty
r^{d-1-\frac{2(d+2s)}{p-1}}\,dr<\infty.
\]
The latter condition holds if and only if $p<3+4s/d$. Therefore, \eqref{eq:L1-threshold} is proved.
\end{rem}

Finally, we establish a strict binding inequality for the optimal
constants in the finite-rank fractional Lieb--Thirring inequality.

\begin{thm}\label{thm:doubling}
Let $d\in\mathbb N$, $0<s<1$, and assume that
\begin{equation}\label{eq:doubling-range}
	\gamma>0,\qquad
	p=\gamma+\frac{d}{2s}\geq3+\frac{4s}{d}.
\end{equation}
Then, for every $N\in\mathbb N$, we have 
\begin{equation}\label{eq:doubling}
	L^{(2N)}_{\gamma,d,s}>L^{(N)}_{\gamma,d,s}.
\end{equation}
In particular, $L^{(N)}_{\gamma,d,s}<L_{\gamma,d,s}$
for every $N\in\mathbb N$.

Moreover, every normalized maximizing sequence for
$L^{(N)}_{\gamma,d,s}$ is precompact modulo translations
in $L^p(\mathbb R^d)$.	
\end{thm}

\begin{rem}
Under the assumptions of Theorem~\ref{thm:doubling},
let $V$ be an optimizer for $L^{(N)}_{\gamma,d,s}$,
and let $K$ be as in Theorem~\ref{thm:structure}.
Then
\[
K>N/2.
\]
Indeed, if $2K\leq N$, then
\eqref{eq:effective-rank}, \eqref{eq:doubling}, and
\eqref{eq:monotonicity} yield
$L^{(N)}_{\gamma,d,s}
=L^{(K)}_{\gamma,d,s}
<L^{(2K)}_{\gamma,d,s}
\leq L^{(N)}_{\gamma,d,s}$,
a contradiction.

Moreover, if the optimal constant $L_{\gamma,d,s}$ were attained
by a potential with finitely many negative eigenvalues,
say $M\geq1$, then
$L_{\gamma,d,s}\leq L^{(M)}_{\gamma,d,s}
<L_{\gamma,d,s}$,
a contradiction. Thus every optimizer for
$L_{\gamma,d,s}$ must have infinitely many negative eigenvalues.
\end{rem}

\begin{rem}
The condition \eqref{eq:doubling-range} is sufficient
for the present argument and is not claimed to be optimal.
Indeed, when $p<3+4s/d$, the estimate
\eqref{eq:overlap-integral} no longer guarantees
$R^{-2(d+2s)}=o(A_R)$, which is needed to control
the remainder in the trace expansion
(see \eqref{eq:matrix-norm-estimate} and
\eqref{eq:numerator}).	
\end{rem}

Although finite-rank optimizers have been studied for the local
Schr\"odinger operator \cite{FGL25}, extending these results to
fractional Schr\"odinger operator introduces several additional
difficulties. The main analytical challenges in the fractional
setting are the following.
\begin{enumerate}
\item The existence of finite-rank optimizers in the fractional setting does not follow directly from the corresponding local theory.
The main difficulty comes from the nonlocality of the fractional Laplacian. Indeed, the associated quadratic forms contain nonvanishing interaction terms.
To overcome this difficulty, we establish a fractional spectral splitting result \eqref{eq:splitting} by combining fractional IMS localization, quantitative estimates of the nonlocal interaction terms, and orthogonalization arguments.
This result allows us to establish the existence of optimizers.

\item The strict binding inequality also requires a different analysis from the local case. 
Due to the algebraic decay of fractional eigenfunctions, the interaction generated by two translated copies has to be compared
carefully with the error terms arising from the spectral trace
expansion and the orthogonalization procedure. 
By combining the Euler--Lagrange equation with resolvent kernel
estimates, we obtain the precise asymptotic behavior of optimizers
at infinity, which allows us to quantify the interaction terms.
Comparing these leading interaction contributions with the error terms enables us to prove the strict inequality between finite-rank constants.
\end{enumerate}

This paper is organized as follows.
Section~\ref{sec:preliminaries} collects some preliminary results and introduces the spectral framework for fractional Schr\"odinger operators.
In Section~\ref{sec:existence}, we prove the existence of optimizers by analyzing the spectrum in the case of dichotomy.
In Section~\ref{sec:structure}, we establish the spectral gap, derive the Euler--Lagrange equation, and investigate the regularity and asymptotic properties of optimizers and the associated eigenfunctions.
Finally, in Section~\ref{sec:doubling}, we derive the strict binding inequality via a trace expansion and obtain precompactness of normalized maximizing sequences modulo translations.

We conclude with the following notation:
\begin{itemize}
	\item For a measurable set $A\subset\mathbb R^d$, $\mathbf{1}_A$ denotes its characteristic function.
	
	\item The symbol $\langle\cdot,\cdot\rangle$ denotes the
	$L^2(\mathbb R^d)$ inner product unless otherwise indicated.

	\item Unless otherwise stated, limits of sequences are taken as their indices tend to infinity.
	
	\item The symbols $\|\cdot\|_p$ and $\|\cdot\|_{H^s}$ denote the usual norms on	$L^p(\mathbb R^d)$ and $H^s(\mathbb R^d)$,	respectively.
	
	\item $B_r:=B_r(0)=\{x\in\mathbb R^d:|x|<r\}$ denotes the open ball centered at the origin with radius $r$.
	
	\item Throughout the paper, $C,C_0,C_1,\ldots$ and $c_0,c_1,\ldots$ denote positive constants whose values may change from line to line unless otherwise specified. Dependence on parameters is indicated when relevant.
\end{itemize}

\section{Preliminaries}\label{sec:preliminaries}
In this section, we collect several preliminary results that will be used repeatedly in the sequel. We first recall the sharp fractional Gagliardo--Nirenberg--Sobolev inequality, its optimizer characterization, and the kernel representation of the fractional kinetic energy; see \cite[Proposition~3.1 and Theorem~3.4]{FLS16}, \cite[Section~2.2]{Din20}, and \cite[Proposition~3.4]{DPV12}.
\begin{lem}\label{lem:GNS}
	Let $d\geq 1$, $s\in(0,1)$, $u\in H^s(\mathbb{R}^d)$, and
	\[
	2<\alpha<2_s^*:=\begin{cases}
		\dfrac{2d}{d-2s}, & d>2s,\\[1mm]
		+\infty, & d\leq 2s.
	\end{cases}
	\]
	Then the sharp fractional Gagliardo--Nirenberg--Sobolev inequality
	\begin{equation*}
		\|u\|_{\alpha}^{\alpha}
		\leq
		C_{\mathrm{opt}}
		\left\|(-\Delta)^{\frac{s}{2}}u\right\|_2^{\frac{d(\alpha-2)}{2s}}
		\|u\|_2^{\alpha-\frac{d(\alpha-2)}{2s}}
	\end{equation*}
	holds with the optimal constant $C_{\mathrm{opt}}=C(s,d,\alpha)$.
	Moreover, $C_{\mathrm{opt}}$ is attained by the unique positive radial
	ground state $Q_{d,\alpha}$ of $(-\Delta)^sQ+Q=|Q|^{\alpha-2}Q$ in $\mathbb{R}^d$.
	In particular,
	\[
	C_{\mathrm{opt}}
	=
	\left(
	\frac{2d-(d-2s)\alpha}{d(\alpha-2)}
	\right)^{\frac{d(\alpha-2)}{4s}}
	\frac{2s\alpha}{2d-(d-2s)\alpha}
	\frac{1}{\|Q_{d,\alpha}\|_{2}^{\alpha-2}}.
	\]
	
	Furthermore, the fractional kinetic energy admits the kernel representation
	\begin{equation}\label{eq:kernel-form}
		\left\|(-\Delta)^{s/2}u\right\|_2^2
		=
		\frac{c_{d,s}}{2}
		\iint_{\mathbb{R}^d\times\mathbb{R}^d}
		\frac{|u(x)-u(y)|^2}{|x-y|^{d+2s}}\,dx\,dy,
	\end{equation}
	where $c_{d,s}>0$ depends only on $d$ and $s$.
\end{lem}

We next record the following fractional Sobolev embeddings; see
\cite[Theorems~6.5, 6.9 and~8.2]{DPV12} and \cite[Lemma~2.1]{Din20}.
\begin{lem}\label{lem:fractional-embedding}
	Let $d\geq 1$ and $0<s<1$. Then the following embeddings are continuous:
	\begin{itemize}
		\item $H^s(\mathbb{R}^d)\hookrightarrow L^q(\mathbb{R}^d)$ for every
		$2\leq q\leq \frac{2d}{d-2s}$ when $d>2s$;
		\item $H^s(\mathbb{R}^d)\hookrightarrow L^q(\mathbb{R}^d)$ for every
		$2\leq q<\infty$ when $d=2s$;
		\item $H^s(\mathbb{R}^d)\hookrightarrow C_b(\mathbb{R}^d)$ when $d<2s$,
		where $C_b(\mathbb{R}^d)$ denotes the space of bounded continuous
		functions on $\mathbb{R}^d$.
	\end{itemize}
	Moreover, the embedding $H^s(\mathbb{R}^d)\hookrightarrow L^q_{\mathrm{loc}}(\mathbb{R}^d)$ is locally compact for every $1\leq q<\frac{2d}{d-2s}$ when $d>2s$, and for every $1\leq q<\infty$ when $d\leq 2s$.	
\end{lem}

Finally, we introduce the operator-theoretic framework for the fractional  Schr\"odinger operators considered in this paper, providing the necessary spectral setting for the analysis of negative eigenvalues.
\begin{lem}\label{lem:form-compactness}
	For $d\geq1$, $0<s<1$, $\gamma>\max\{0,1-\frac{d}{2s}\}$, and
	$p=\gamma+\frac{d}{2s}$, assume that $V\in L^p(\mathbb R^d)$ is nonnegative. Then the quadratic form
	\[
	q_V[u]=\|(-\Delta)^{s/2}u\|_2^2-\int_{\mathbb R^d}V|u|^2\,dx,
	\qquad u\in H^s(\mathbb R^d),
	\]
	defines a self-adjoint operator $H_V=(-\Delta)^s-V$ with $\sigma_{\mathrm{ess}}(H_V)=[0,\infty)$.	
\end{lem}
\begin{proof}
	Let
	\[
	p':=\frac{p}{p-1},
	\qquad
	\vartheta:=\frac{d}{2sp}
	=\frac{d}{d+2s\gamma}\in(0,1).
	\]
	The assumption on $\gamma$ ensures that $p>1$ and $2<2p'<2_s^*$. Hence, Lemma~\ref{lem:GNS} with $\alpha=2p'$ gives	
	\begin{equation}\label{eq:GN}
		\|u\|_{2p'}^2
		\leq
		C
		\left\|(-\Delta)^{s/2}u\right\|_2^{2\vartheta}
		\|u\|_2^{2(1-\vartheta)}.
	\end{equation}
	For $u\in H^s(\mathbb{R}^d)$ and $\varepsilon>0$, it follows from
	H\"older's inequality, \eqref{eq:GN}, and Young's inequality that
	\begin{equation}\label{eq:form}
		\begin{aligned}
			\int_{\mathbb{R}^d}V|u|^2\,dx
			\leq
			\|V\|_p\|u\|_{2p'}^2 
			&\leq
			C\|V\|_p
			\left\|(-\Delta)^{s/2}u\right\|_2^{2\vartheta}
			\|u\|_2^{2(1-\vartheta)} \\
			&\leq
			\varepsilon\left\|(-\Delta)^{s/2}u\right\|_2^2
			+
			C\varepsilon^{-d/(2s\gamma)}
			\|V\|_p^{p/\gamma}
			\|u\|_2^2.
		\end{aligned}
	\end{equation}
	Thus, $V$ is infinitesimally form bounded with respect to $(-\Delta)^s$. By the KLMN theorem \cite[Theorem X.17]{RSII75}, the associated quadratic form defines a self-adjoint operator $H_V$. 
	
	We next show that $V$ is form compact with respect to $(-\Delta)^s$. To see this, choose compactly supported $V_k\in L^\infty(\mathbb{R}^d)$ such that $\|V_k-V\|_p\to0$. For fixed $k$, let $\Omega\supset\operatorname{supp}V_k$ be a bounded open set.
	If $u_n\rightharpoonup0$ in $H^s(\mathbb{R}^d)$, then Lemma~\ref{lem:fractional-embedding} implies $u_n\to0$ in $L^2(\Omega)$. Consequently,
	\begin{equation}\label{eq2.2}
		\left|
		\int_{\mathbb{R}^d}V_k|u_n|^2\,dx
		\right|
		\leq
		\|V_k\|_\infty
		\|u_n\|_{L^2(\Omega)}^2
		\to0,
	\end{equation}
	so $V_k$ is form compact with respect to $(-\Delta)^s$. Furthermore, H\"older's inequality and \eqref{eq:GN} imply
	\begin{equation}\label{eq2.3}
		\left|
		\int_{\mathbb{R}^d}(V-V_k)|u_n|^2\,dx
		\right|
		\leq
		\|V-V_k\|_p\|u_n\|_{2p'}^2
		\leq
		C\|V-V_k\|_p
		\left\|(-\Delta)^{s/2}u_n\right\|_2^{2\vartheta}
		\|u_n\|_2^{2(1-\vartheta)}.
	\end{equation}
	Since $\{u_n\}$ is bounded in $H^s(\mathbb{R}^d)$, by \eqref{eq2.2} and \eqref{eq2.3}, we get
	\[
	\left|
	\int_{\mathbb{R}^d}V|u_n|^2\,dx
	\right|
	\leq
	\left|
	\int_{\mathbb{R}^d}V_k|u_n|^2\,dx
	\right|
	+
	C\|V-V_k\|_p.
	\]
	Letting first $n\to\infty$ and then $k\to\infty$, we obtain $\int_{\mathbb{R}^d}V|u_n|^2\,dx\to0$, and the desired form compactness follows. By Weyl's theorem \cite[Theorem XIII.14]{RSIV78}, one can see that $\sigma_{\mathrm{ess}}(H_V)=
	\sigma_{\mathrm{ess}}((-\Delta)^s)=[0,\infty)$.
\end{proof}

\section{Existence of optimizers for the finite-rank problem}\label{sec:existence}
In this section, we aim to prove Theorem~\ref{thm:existence}, namely the existence of optimizers for the finite-rank Lieb--Thirring constants $L^{(N)}_{\gamma,d,s}$. 
In Section~\ref{subsec:forms}, we present several spectral and localization estimates, as well as compactness results that will be used below.
In Section~\ref{subsec:splitting}, we analyze the spectrum
in the case of dichotomy.
Finally, in Section~\ref{subsec:compactness}, we apply the concentration--compactness argument to normalized maximizing sequences to exclude vanishing, treat the possible dichotomy, and thereby prove the existence of optimizers.

\subsection{Spectral estimates, localization, and compactness} \label{subsec:forms}
We begin with some basic spectral estimates.

\begin{lem}\label{lem:spectral-stability}
	Let $V\geq0$ belong to $L^p(\mathbb{R}^d)$. Then
	\begin{equation}\label{eq:ground-bound}
		|\lambda_1(H_V)|^\gamma
		\leq
		C\|V\|_p^p.
	\end{equation}
	
	Moreover, let $U\geq0$ belong to $L^p(\mathbb{R}^d)$ and $j\geq1$.
	For every $M>0$ with $\|U\|_p,\|V\|_p\leq M$, we have
	\begin{equation}\label{eq:level-Lipschitz}
		|\lambda_j(H_U)-\lambda_j(H_V)|
		\leq
		C_M\|U-V\|_p,
	\end{equation}
	and, uniformly for $1/2\leq c\leq1$,
	\begin{equation}\label{eq:kinetic-continuity}
		0\leq
		\lambda_j((-\Delta)^s-V)
		-\lambda_j(c(-\Delta)^s-V)
		\leq
		C_M(1-c).
	\end{equation}
\end{lem}

\begin{proof}
	Taking $\varepsilon=1/2$ in \eqref{eq:form}, one has
	\[
	q_V[u]
	\geq
	\frac12\left\|(-\Delta)^{s/2}u\right\|_2^2
	-C\|V\|_p^{p/\gamma}\|u\|_2^2
	\geq
	-C\|V\|_p^{p/\gamma}\|u\|_2^2.
	\]
	Hence, by the min--max principle,
	\[
	-C\|V\|_p^{p/\gamma}
	\leq
	\lambda_1(H_V)
	\leq 0,
	\]
	which proves \eqref{eq:ground-bound}.
	
	To prove \eqref{eq:level-Lipschitz}, we first establish a uniform kinetic bound. For $W\in\{U,V\}$ and
	$a\in[1/2,1]$, we introduce
	\[
	H_{W,a}:=a(-\Delta)^s-W,
	\qquad
	q_{a,W}[u]
	:=
	a\left\|(-\Delta)^{s/2}u\right\|_2^2
	-\int_{\mathbb{R}^d}W|u|^2\,dx.
	\]
	If $u$ belongs to the negative spectral subspace of $H_{W,a}$, then
	$q_{a,W}[u]\leq0$. 
	Setting $\varepsilon=1/4$ in \eqref{eq:form}, we obtain
	\[
	\begin{aligned}
		a\left\|(-\Delta)^{s/2}u\right\|_2^2
		\leq
		\int_{\mathbb{R}^d}W|u|^2\,dx 
		\leq
		\frac14\left\|(-\Delta)^{s/2}u\right\|_2^2
		+
		C\|W\|_p^{p/\gamma}\|u\|_2^2 
		\leq
		\frac14\left\|(-\Delta)^{s/2}u\right\|_2^2
		+
		C_M\|u\|_2^2,
	\end{aligned}
	\]
	which, since $a\geq1/2$, implies
	\begin{equation}\label{eq:uniform-kinetic}
		\left\|(-\Delta)^{s/2}u\right\|_2^2
		\leq
		C_M\|u\|_2^2.
	\end{equation}
	On every negative spectral subspace of $H_U$ or $H_V$, we infer from H\"older's inequality, \eqref{eq:GN}, and \eqref{eq:uniform-kinetic} that
	\[
	\begin{aligned}
		|q_U[u]-q_V[u]|
		\leq
		\|U-V\|_p\|u\|_{2p'}^2 
		\leq
		C\|U-V\|_p
		\left\|(-\Delta)^{s/2}u\right\|_2^{2\vartheta}
		\|u\|_2^{2(1-\vartheta)} 
		\leq
		C_M\|U-V\|_p\|u\|_2^2.
	\end{aligned}
	\]
	If $\lambda_j(H_U)=\lambda_j(H_V)=0$, then \eqref{eq:level-Lipschitz} is immediate. Otherwise, assume without loss of generality that $\lambda_j(H_U)<0$, and let
	\[
	E_j(U):=\operatorname{span}\{\phi_1^U,\ldots,\phi_j^U\},
	\]
	where $\phi_1^U,\ldots,\phi_j^U$ are orthonormal eigenfunctions
	corresponding to the first $j$ negative eigenvalues of $H_U$.
	Then, by the min--max principle and the preceding estimate,
	\[
	\begin{aligned}
		\lambda_j(H_V)
		\leq
		\max_{\substack{u\in E_j(U)\\ \|u\|_2=1}}q_V[u] 
		\leq
		\max_{\substack{u\in E_j(U)\\ \|u\|_2=1}}
		\left(q_U[u]+C_M\|U-V\|_p\right) 
		=
		\lambda_j(H_U)+C_M\|U-V\|_p.
	\end{aligned}
	\]
	If $\lambda_j(H_V)<0$, interchanging $U$ and $V$ yields the reverse inequality. If $\lambda_j(H_V)=0$, the above estimate gives $-\lambda_j(H_U)\leq C_M\|U-V\|_p$. Thus \eqref{eq:level-Lipschitz} follows.
	
	Finally, it remains to establish \eqref{eq:kinetic-continuity}. Fix $c\in[1/2,1]$. If $\lambda_j(H_{U,c})<0$, let
	$E_j(c)$ denote the span of the first $j$ negative eigenfunctions of
	$H_{U,c}$. Since
	\[
	q_{1,U}[u]
	=
	q_{c,U}[u]
	+
	(1-c)\left\|(-\Delta)^{s/2}u\right\|_2^2,
	\]
	it follows from \eqref{eq:uniform-kinetic} and the min--max principle that
	\[
	\begin{aligned}
		\lambda_j(H_{U,1})
		\leq
		\max_{\substack{u\in E_j(c)\\ \|u\|_2=1}}q_{1,U}[u] 
		&=
		\max_{\substack{u\in E_j(c)\\ \|u\|_2=1}}
		\left(
		q_{c,U}[u]
		+
		(1-c)\left\|(-\Delta)^{s/2}u\right\|_2^2
		\right) \\
		&\leq
		\max_{\substack{u\in E_j(c)\\ \|u\|_2=1}}
		q_{c,U}[u]
		+
		C_M(1-c) \\
		&=
		\lambda_j(H_{U,c})+C_M(1-c).
	\end{aligned}
	\]
	On the other hand, since $q_{c,U}[u]\leq q_{1,U}[u]$ for all
	$u\in H^s(\mathbb{R}^d)$, the min--max principle gives
	\begin{equation}\label{eq:kinetic-monotonicity}
		\lambda_j(H_{U,c})\leq \lambda_j(H_{U,1}).
	\end{equation}
	This proves \eqref{eq:kinetic-continuity} when
	$\lambda_j(H_{U,c})<0$.
	Moreover, if $\lambda_j(H_{U,c})=0$, then
	\eqref{eq:kinetic-monotonicity} implies
	$\lambda_j(H_{U,1})=0$, and the conclusion follows directly. Consequently,
	\eqref{eq:kinetic-continuity} holds in all cases.
\end{proof}

We will use the following discrete and continuous IMS-type localization estimates for the fractional kinetic energy.
\begin{lem}\label{lem:fractional-localization}
	For $R>0$ and $u\in H^s(\mathbb{R}^d)$, the following statements hold.
	\begin{enumerate}[label=(\roman*),ref=(\roman*)]
		\item\label{item:discrete-localization}
		Let $\{\chi_\ell\}_{\ell=1}^L$ be a finite family of real-valued smooth functions on $\mathbb{R}^d$ satisfying
		\[
		\sum_{\ell=1}^L\chi_\ell^2=1,
		\qquad
		\sum_{\ell=1}^L|\chi_\ell(x)-\chi_\ell(y)|^2
		\leq
		C_0\min\left\{1,\frac{|x-y|^2}{R^2}\right\},
		\]
		where $C_0>0$ is independent of $R$. Then there exists $C>0$, independent of $R$, such that
		\begin{equation}\label{eq:IMS}
			\left|
			\sum_{\ell=1}^L
			\left\|(-\Delta)^{s/2}(\chi_\ell u)\right\|_2^2
			-
			\left\|(-\Delta)^{s/2}u\right\|_2^2
			\right|
			\leq
			CR^{-2s}\|u\|_2^2.
		\end{equation}
		
		\item\label{item:continuous-localization}
		Let $\chi\in C_c^\infty(B_1)$ be real-valued and satisfy 
		$\|\chi\|_2^2=1$. For $z\in\mathbb{R}^d$, set
		\begin{equation}\label{eq:chi-Rz}
			\chi_{R,z}(x)
			:=
			R^{-d/2}\chi\left(\frac{x-z}{R}\right).
		\end{equation}
		Then, for every $x\in\mathbb{R}^d$,
		$\int_{\mathbb{R}^d}\chi_{R,z}^2(x)\,dz=1$, and
		\begin{equation}\label{eq:continuous-IMS}
			\left|
			\int_{\mathbb{R}^d}
			\left\|(-\Delta)^{s/2}(\chi_{R,z}u)\right\|_2^2\,dz
			-
			\left\|(-\Delta)^{s/2}u\right\|_2^2
			\right|
			\leq
			CR^{-2s}\|u\|_2^2,
		\end{equation}
		where $C>0$ is independent of $R$.
	\end{enumerate}	
\end{lem}

\begin{proof}
	Since $\chi_\ell$ is real-valued and $\sum_\ell^L\chi_\ell^2=1$, it follows that
	\[
	\sum_\ell^L|\chi_\ell(x)u(x)-\chi_\ell(y)u(y)|^2
	=
	|u(x)|^2+|u(y)|^2
	-2\operatorname{Re}\bigl(u(x)\overline{u(y)}\bigr)
	\sum_\ell^L\chi_\ell(x)\chi_\ell(y).
	\]
	We also note that
	\[
	\sum_\ell^L|\chi_\ell(x)-\chi_\ell(y)|^2
	=
	2-2\sum_\ell^L\chi_\ell(x)\chi_\ell(y).
	\]
	Together with $|u(x)-u(y)|^2
	=
	|u(x)|^2+|u(y)|^2
	-2\operatorname{Re}\bigl(u(x)\overline{u(y)}\bigr)$, these identities yield
	\[
	\sum_\ell^L|\chi_\ell(x)u(x)-\chi_\ell(y)u(y)|^2-|u(x)-u(y)|^2
	=
	\operatorname{Re}\bigl(u(x)\overline{u(y)}\bigr)
	\sum_\ell^L|\chi_\ell(x)-\chi_\ell(y)|^2.
	\]
	Thus, by \eqref{eq:kernel-form}, we infer that
	\begin{equation}\label{eq:IMS-kernel-identity}
		\sum_\ell^L\|(-\Delta)^{s/2}(\chi_\ell u)\|_2^2
		-\|(-\Delta)^{s/2}u\|_2^2 =
		\frac{c_{d,s}}{2}
		\iint_{\mathbb{R}^d\times\mathbb{R}^d}
		\frac{\operatorname{Re}\bigl(u(x)\overline{u(y)}\bigr)}
		{|x-y|^{d+2s}}
		\sum_\ell^L|\chi_\ell(x)-\chi_\ell(y)|^2\,dx\,dy,
	\end{equation}
	which, combined with the fact that
	$\sum_\ell^L|\chi_\ell(x)-\chi_\ell(y)|^2
	\leq C\min\left\{1,\frac{|x-y|^2}{R^2}\right\}$,
	means that
	\[
	\begin{aligned}
		\left|
		\sum_\ell^L\|(-\Delta)^{s/2}(\chi_\ell u)\|_2^2
		-\|(-\Delta)^{s/2}u\|_2^2
		\right|
		&\leq
		C\iint_{\mathbb{R}^d\times\mathbb{R}^d}
		\frac{|u(x)u(y)|}{|x-y|^{d+2s}}
		\min\left\{1,\frac{|x-y|^2}{R^2}\right\}\,dx\,dy.
	\end{aligned}
	\]
	Since $2|u(x)u(y)|\leq |u(x)|^2+|u(y)|^2$, using the symmetry of the kernel and Tonelli's theorem, we obtain
	\[
	\begin{aligned}
		&\iint_{\mathbb{R}^d\times\mathbb{R}^d}
		\frac{|u(x)u(y)|}{|x-y|^{d+2s}}
		\min\left\{1,\frac{|x-y|^2}{R^2}\right\}\,dx\,dy \\
		&\leq
		\frac12
		\iint_{\mathbb{R}^d\times\mathbb{R}^d}
		\frac{|u(x)|^2+|u(y)|^2}{|x-y|^{d+2s}}
		\min\left\{1,\frac{|x-y|^2}{R^2}\right\}\,dx\,dy \\
		&=
		\iint_{\mathbb{R}^d\times\mathbb{R}^d}
		\frac{|u(x)|^2}{|x-y|^{d+2s}}
		\min\left\{1,\frac{|x-y|^2}{R^2}\right\}\,dx\,dy 
		=\|u\|_2^2
		\int_{\mathbb{R}^d}
		\frac{\min\{1,|z|^2/R^2\}}{|z|^{d+2s}}\,dz.
	\end{aligned}
	\]
	Furthermore, setting $z:=Rw$, one finds
	\begin{equation}\label{eq:IMS-kernel-scaling}
		\begin{aligned}
			\int_{\mathbb{R}^d}\frac{\min\{1,|z|^2/R^2\}}{|z|^{d+2s}}\,dz
			&=
			R^{-2s}\int_{\mathbb{R}^d}\frac{\min\{1,|w|^2\}}{|w|^{d+2s}}\,dw \\
			&\le C R^{-2s}
			\left(\int_0^1 r^{1-2s}\,dr+\int_1^\infty r^{-1-2s}\,dr\right)
			\leq CR^{-2s},
		\end{aligned}
	\end{equation}
	where the last two integrals are finite since $0<s<1$. Therefore, \eqref{eq:IMS} is established, completing the proof of \ref{item:discrete-localization}.
	
	We now prove \ref{item:continuous-localization}. By the normalization of $\chi$, it follows immediately that
	\begin{equation}\label{eq:chi-normalization}
		\int_{\mathbb{R}^d}\chi_{R,z}^2(x)\,dz
		=
		\int_{\mathbb{R}^d}\chi^2(w)\,dw
		=
		1.
	\end{equation}
	First, to control large values of $|x-y|$, after the change of variables $w:=(x-z)/R$, we get
	\begin{equation}\label{eq:large-distance}
		\int_{\mathbb{R}^d}|\chi_{R,z}(x)-\chi_{R,z}(y)|^2\,dz
		=
		\int_{\mathbb{R}^d}
		\left|
		\chi(w)-\chi\left(w+\frac{y-x}{R}\right)
		\right|^2\,dw
		\leq 4.
	\end{equation}
	On the other hand, to control small values of $|x-y|$, the fundamental theorem of calculus gives
	\[
	\chi(w+h)-\chi(w)
	=
	\int_0^1 h\cdot\nabla\chi(w+th)\,dt,
	\qquad
	h=\frac{y-x}{R},
	\]
	which, together with the Cauchy--Schwarz inequality, implies that
	\begin{equation}\label{eq:small-distance}
		\begin{aligned}
			\int_{\mathbb{R}^d}|\chi(w+h)-\chi(w)|^2\,dw
			&\leq
			|h|^2\int_{\mathbb{R}^d}\int_0^1
			|\nabla\chi(w+th)|^2\,dt\,dw \\
			&=
			|h|^2\int_0^1\|\nabla\chi\|_2^2\,dt 
			=
			\frac{|x-y|^2}{R^2}\|\nabla\chi\|_2^2 
			\leq
			C\frac{|x-y|^2}{R^2}.
		\end{aligned}
	\end{equation}
	Consequently, from \eqref{eq:large-distance} and
	\eqref{eq:small-distance}, we deduce that $\int_{\mathbb{R}^d}|\chi_{R,z}(x)-\chi_{R,z}(y)|^2\,dz
	\leq
	C\min\left\{1,\frac{|x-y|^2}{R^2}\right\}.$
	Using \eqref{eq:chi-normalization} and arguing as in \ref{item:discrete-localization}, we obtain
	\[
	\int_{\mathbb{R}^d}
	|\chi_{R,z}(x)u(x)-\chi_{R,z}(y)u(y)|^2\,dz
	-|u(x)-u(y)|^2=
	\operatorname{Re}\bigl(u(x)\overline{u(y)}\bigr)
	\int_{\mathbb{R}^d}
	|\chi_{R,z}(x)-\chi_{R,z}(y)|^2\,dz.
	\]
	Together with the preceding estimate and \eqref{eq:kernel-form}, this yields
	\[
	\left|
	\int_{\mathbb{R}^d}
	\|(-\Delta)^{s/2}(\chi_{R,z}u)\|_2^2\,dz
	-\|(-\Delta)^{s/2}u\|_2^2
	\right|
	\leq CR^{-2s}\|u\|_2^2,
	\]
	thereby proving \eqref{eq:continuous-IMS} and concluding the proof.
\end{proof}

We next establish a nonvanishing criterion and a convergence result for the min--max levels.
\begin{lem}\label{lem:tight}
	The following two assertions hold.
	\begin{enumerate}[label=(\roman*),ref=(\roman*)]
		\item\label{item:nonvanishing}
		If $V\geq0$ and $\|V\|_p^p\leq1$, then
		\begin{equation}\label{eq:nonvanishing}
			|\lambda_1(H_V)|
			\leq
			C\left(
			\sup_{y\in\mathbb{R}^d}
			\int_{B_1(y)}V^p\,dx
			\right)^{1/p}.
		\end{equation}
		
		\item\label{item:eigenvalue-convergence}
		Suppose that $V_n\geq0$, $V_n\rightharpoonup V$ in
		$L^p(\mathbb{R}^d)$, and
		\begin{equation}\label{eq:tight}
			\lim_{R\to\infty}\limsup_{n\to\infty}
			\int_{|x|>R}V_n^p\,dx=0.
		\end{equation}
		Then $\lambda_j(H_{V_n})\to\lambda_j(H_V)$ for every fixed $j\geq1$.
	\end{enumerate}	
\end{lem}

\begin{proof}
	We start by proving \ref{item:nonvanishing}. Define
	\[
	m:=\sup_{y\in\mathbb{R}^d}\int_{B_1(y)}V^p\,dx\in[0,1].
	\]
	For $R\geq1$, the ball $B_R(z)$ can be covered by at most $CR^d$
	unit balls, and hence $\int_{B_R(z)}V^p\,dx\leq CR^dm$. Let $\chi_{R,z}$ be defined by \eqref{eq:chi-Rz}. Since
	$\operatorname{supp}\chi_{R,z}\subset B_R(z)$, applying
	\eqref{eq:form} to $V\mathbf{1}_{B_R(z)}$ and $\chi_{R,z}u$ with
	$\varepsilon=1$, we obtain
	\begin{equation}\label{eq:localized-form-bound}
		\int_{\mathbb{R}^d}V\chi_{R,z}^2|u|^2\,dx
		\leq
		\left\|(-\Delta)^{s/2}(\chi_{R,z}u)\right\|_2^2
		+
		CR^{d/\gamma}m^{1/\gamma}\|\chi_{R,z}u\|_2^2.
	\end{equation}
	From \eqref{eq:chi-normalization} and Tonelli's theorem, it follows that
	\[
	\int_{\mathbb{R}^d}\int_{\mathbb{R}^d}
	V(x)\chi_{R,z}^2(x)|u(x)|^2\,dx\,dz
	=
	\int_{\mathbb{R}^d}V|u|^2\,dx,
	\qquad
	\int_{\mathbb{R}^d}\|\chi_{R,z}u\|_2^2\,dz
	=
	\|u\|_2^2.
	\]
	Integrating \eqref{eq:localized-form-bound} with respect to $z$ and
	combining these identities with \eqref{eq:continuous-IMS}, we deduce that
	\[
	\begin{aligned}
		\int_{\mathbb{R}^d}V|u|^2\,dx
		&\leq
		\int_{\mathbb{R}^d}
		\left\|(-\Delta)^{s/2}(\chi_{R,z}u)\right\|_2^2\,dz
		+
		CR^{d/\gamma}m^{1/\gamma}\|u\|_2^2 \\
		&\leq
		\left\|(-\Delta)^{s/2}u\right\|_2^2
		+
		CR^{-2s}\|u\|_2^2
		+
		CR^{d/\gamma}m^{1/\gamma}\|u\|_2^2 \\
		&\leq
		\left\|(-\Delta)^{s/2}u\right\|_2^2
		+
		C\left(R^{-2s}+R^{d/\gamma}m^{1/\gamma}\right)\|u\|_2^2.
	\end{aligned}
	\]
	Therefore,
	\[
	q_V[u]
	\geq
	-C\left(R^{-2s}+R^{d/\gamma}m^{1/\gamma}\right)\|u\|_2^2.
	\]
	Suppose first that $m>0$. Since $m\leq1$, the choice
	$R=m^{-1/(2sp)}$ satisfies $R\geq1$. Recalling that
	$p=\gamma+d/(2s)$, we obtain
	$R^{-2s}=m^{1/p}$ and
	$R^{d/\gamma}m^{1/\gamma}
	=m^{1/\gamma-d/(2sp\gamma)}
	=m^{1/p}$.
	Hence
	\[
	q_V[u]\geq-Cm^{1/p}\|u\|_2^2.
	\]
	By the min--max principle,
	$-Cm^{1/p}\leq\lambda_1(H_V)\leq0$, and therefore
	\[
	|\lambda_1(H_V)|
	\leq
	Cm^{1/p}
	=
	C\left(
	\sup_{y\in\mathbb{R}^d}
	\int_{B_1(y)}V^p\,dx
	\right)^{1/p}.
	\]
	If $m=0$, then $\int_{B_1(y)}V^p\,dx=0$ for every $y\in\mathbb{R}^d$, so $V=0$ almost everywhere in $\mathbb{R}^d$, and hence $\lambda_1(H_V)=0$. Consequently, \eqref{eq:nonvanishing} holds in all cases.
	
	We turn to \ref{item:eigenvalue-convergence} and claim that
	\begin{equation}\label{eq:weak-form-norm}
		\sup_{\|u\|_{H^s}\leq1}
		\left|
		\int_{\mathbb{R}^d}(V_n-V)|u|^2\,dx
		\right|\to 0.
	\end{equation} 
	Indeed, by Lemma~\ref{lem:fractional-embedding}, the embedding
	$H^s(\mathbb{R}^d)\hookrightarrow L^{2p'}(B_R)$ is compact for every
	fixed $R>0$, where $p'=\frac{p}{p-1}$. It then follows that
	\[
	\mathcal K_R:=
	\left\{|u|^2:\ u\in H^s(\mathbb{R}^d),\ \|u\|_{H^s}\leq1\right\}
	\]
	is relatively compact in $L^{p'}(B_R)$. Given $\varepsilon>0$, choose
	$f_1,\ldots,f_K\in\mathcal K_R$ such that, for every
	$f\in\mathcal K_R$, there exists $k\in\{1,\ldots,K\}$ with $\|f-f_k\|_{L^{p'}(B_R)}<\varepsilon$. Moreover, since
	$V_n\rightharpoonup V$ in $L^p(\mathbb{R}^d)$, the sequence
	$\{V_n-V\}$ is bounded in $L^p(\mathbb{R}^d)$. Hence
	\[
	\begin{aligned}
		\left|\int_{B_R}(V_n-V)f\,dx\right|
	&\leq
	\left|\int_{B_R}(V_n-V)f_k\,dx\right|
	+
	\|V_n-V\|_p\|f-f_k\|_{L^{p'}(B_R)}\\
	&\leq
	\left|\int_{B_R}(V_n-V)f_k\,dx\right|
	+C\varepsilon,
	\end{aligned}
	\]
	which, together with $K<\infty$ and $\int_{B_R}(V_n-V)f_k\,dx\to0$ for every $k$, yields 
	\[
	\limsup_{n\to\infty}
	\sup_{f\in\mathcal K_R}
	\left|\int_{B_R}(V_n-V)f\,dx\right|
	\leq C\varepsilon.
	\]
	Passing to the limit $\varepsilon\to0$, we obtain $\sup_{f\in\mathcal K_R}
	\left|\int_{B_R}(V_n-V)f\,dx\right|\to0$. Equivalently,
	\[
	\sup_{\|u\|_{H^s}\leq1}
	\left|
	\int_{B_R}(V_n-V)|u|^2\,dx
	\right|
	\to0.
	\]
	It remains to control the tail. By H\"older's inequality and
	Lemma~\ref{lem:fractional-embedding}, we have
	\[
	\sup_{\|u\|_{H^s}\leq1}
	\left|
	\int_{B_R^c}(V_n-V)|u|^2\,dx
	\right|
	\leq
	C\left(
	\|V_n\mathbf{1}_{B_R^c}\|_p
	+
	\|V\mathbf{1}_{B_R^c}\|_p
	\right).
	\]
	Moreover, weak lower semicontinuity gives $\|V\mathbf{1}_{B_R^c}\|_p
	\leq
	\liminf_{n\to\infty}
	\|V_n\mathbf{1}_{B_R^c}\|_p$. 
	Therefore, we infer from \eqref{eq:tight} that
	\[
	\lim_{R\to\infty}\limsup_{n\to\infty}
	\sup_{\|u\|_{H^s}\leq1}
	\left|
	\int_{B_R^c}(V_n-V)|u|^2\,dx
	\right|
	=0.
	\]
	The local and tail estimates together establish \eqref{eq:weak-form-norm}.
	Let
	\[
	\delta_n:=
	\sup_{\|u\|_{H^s}\leq1}
	\left|
	\int_{\mathbb{R}^d}(V_n-V)|u|^2\,dx
	\right|.
	\]
	By \eqref{eq:weak-form-norm}, $\delta_n\to0$. 
	Moreover, since
	$\sup_n\|V_n\|_p+\|V\|_p<\infty$, \eqref{eq:uniform-kinetic} with
	$a=1$ implies that, on the negative spectral subspaces of $H_{V_n}$
	and $H_V$,
	\[
	|q_{V_n}[u]-q_V[u]|
	=
	\left|
	\int_{\mathbb{R}^d}(V_n-V)|u|^2\,dx
	\right|
	\leq
	\delta_n\|u\|_{H^s}^2
	\leq
	C\delta_n\|u\|_2^2.
	\]
	Using the same min--max argument as in the proof of
	\eqref{eq:level-Lipschitz}, we get
	\[
	|\lambda_j(H_{V_n})-\lambda_j(H_V)|
	\leq C\delta_n\to0
	\]
	for every fixed $j$.
\end{proof}

As a consequence, we have the following compactness result for tight maximizing sequences.
\begin{cor}\label{cor:tight-maximizing}
	Fix $k\in\mathbb{N}$. Let $V_n\geq0$ satisfy \eqref{eq:tight} and
	\begin{equation}\label{eq:tight-maximizing}
		\|V_n\|_p^p\to\alpha>0,
		\qquad
		\sum_{j=1}^k|\lambda_j(H_{V_n})|^\gamma
		\to L^{(k)}_{\gamma,d,s}\alpha.
	\end{equation}
	Then, up to a subsequence, there exists $V\geq0$ such that
	\[
	V_n\to V \qquad\text{in }L^p(\mathbb{R}^d),
	\]
	and $V$ attains the optimal constant
	$L_{\gamma,d,s}^{(k)}$ defined in \eqref{eq:constant}.
\end{cor}

\begin{proof}
	Since $\{V_n\}$ is bounded in $L^p(\mathbb{R}^d)$, after passing to a subsequence, $V_n\rightharpoonup V$ in $L^p(\mathbb{R}^d)$. By Lemma~\ref{lem:tight}-\ref{item:eigenvalue-convergence},
	\eqref{eq:tight-maximizing}, the definition of
	$L^{(k)}_{\gamma,d,s}$, and weak lower semicontinuity, we obtain
	\begin{equation*}
		\begin{aligned}
			L^{(k)}_{\gamma,d,s}\alpha
			=
			\sum_{j=1}^k|\lambda_j(H_{V_n})|^\gamma+o_n(1)
			&=
			\sum_{j=1}^k|\lambda_j(H_V)|^\gamma+o_n(1)\\
			&\leq
			L^{(k)}_{\gamma,d,s}\|V\|_p^p+o_n(1)
			\leq
			L^{(k)}_{\gamma,d,s}\alpha+o_n(1).
		\end{aligned}
	\end{equation*}
	Letting $n\to\infty$, we have $\|V\|_p^p=\alpha$, and hence $V_n\to V$ strongly in $L^p(\mathbb{R}^d)$. Moreover, $V$ is an optimizer for $L^{(k)}_{\gamma,d,s}$.
\end{proof}

\subsection{Spectrum in the case of dichotomy}\label{subsec:splitting}
In this subsection, we establish the following spectral splitting result, which will be used in the concentration--compactness argument.
\begin{lem}\label{lem:splitting}
	Fix $N\in\mathbb{N}$. Let $V_n\geq0$ be bounded in
	$L^p(\mathbb{R}^d)$. Assume that $R_n\to\infty$ and $	\int_{B_{2R_n}\setminus B_{R_n}}V_n^p\,dx\to0$. Set
	\[
	V_n^{\mathrm{in}}:=V_n\mathbf{1}_{B_{R_n}},
	\qquad
	V_n^{\mathrm{out}}:=V_n\mathbf{1}_{\mathbb{R}^d\setminus B_{2R_n}}.
	\]
	Then, up to a subsequence, there exists $k\in\{0,\ldots,N\}$ such that
	\begin{equation}\label{eq:splitting}
		\sum_{j=1}^N|\lambda_j(H_{V_n})|^\gamma
		=
		\sum_{j=1}^k|\lambda_j(H_{V_n^{\mathrm{in}}})|^\gamma
		+
		\sum_{j=1}^{N-k}|\lambda_j(H_{V_n^{\mathrm{out}}})|^\gamma
		+o_n(1).
	\end{equation}	
\end{lem}

\begin{proof}
	The proof is divided into the following five steps.
	
	\textbf{Step 1. Localization and lower spectral estimate.}
	We set
	\begin{equation}\label{eq:splitting-annulus-notation}
		A_n:=B_{2R_n}\setminus B_{R_n},
		\qquad
		\varepsilon_n:=C\|V_n\mathbf{1}_{A_n}\|_p,
	\end{equation}
	where $C>0$ is a fixed constant, independent of $n$, chosen to absorb
	the constants in the estimate below. Clearly, $\varepsilon_n\to0$ as
	$n\to\infty$. Discarding finitely many terms, we may assume that
	$R_n\geq1$ and $\varepsilon_n\leq1/2$. For every
	$u\in H^s(\mathbb{R}^d)$, it follows from H\"older's inequality,
	\eqref{eq:GN}, \eqref{eq:splitting-annulus-notation}, and Young's
	inequality that
	\begin{equation}\label{eq:splitting-annulus-bound}
		\begin{aligned}
			\int_{A_n}V_n|u|^2\,dx
			\leq
			\|V_n\mathbf{1}_{A_n}\|_p\|u\|_{2p'}^2
			&\leq
			C\|V_n\mathbf{1}_{A_n}\|_p
			\|(-\Delta)^{s/2}u\|_2^{2\vartheta}
			\|u\|_2^{2(1-\vartheta)}\\
			&\leq
			\varepsilon_n
			\left(
			\|(-\Delta)^{s/2}u\|_2^2+\|u\|_2^2
			\right).
		\end{aligned}
	\end{equation}
	In addition, choose $\chi_n,\zeta_n,\eta_n\in C^\infty(\mathbb{R}^d;[0,1])$ such that $\chi_n^2+\zeta_n^2+\eta_n^2=1$ and
	\[
	\left\{
	\begin{aligned}
		&\chi_n=1 \quad \text{on } B_{R_n},\\
		&\chi_n=0 \quad \text{on } B_{5R_n/4}^c,\\
		&\|\nabla\chi_n\|_\infty\leq C/R_n,
	\end{aligned}
	\right.
	~~\text{ and }~~
	\left\{
	\begin{aligned}
		&\zeta_n=0 \quad \text{on } B_{7R_n/4},\\
		&\zeta_n=1 \quad \text{on } B_{2R_n}^c,\\
		&\|\nabla\zeta_n\|_\infty+\|\nabla\eta_n\|_\infty\leq C/R_n.
	\end{aligned}
	\right.
	\]
	By the above bounds and the mean value theorem, these cutoffs satisfy the
	hypotheses of \eqref{eq:IMS}. Hence, \eqref{eq:IMS} gives
	\[
	\|(-\Delta)^{s/2}u\|_2^2
	\geq
	\sum_{h\in\{\chi_n,\zeta_n,\eta_n\}}
	\|(-\Delta)^{s/2}(hu)\|_2^2
	-CR_n^{-2s}\|u\|_2^2.
	\]
	Moreover, the support properties imply
	\[
	\int_{\mathbb{R}^d}V_n|u|^2\,dx
	=
	\int_{\mathbb{R}^d}V_n^{\mathrm{in}}|\chi_nu|^2\,dx
	+
	\int_{\mathbb{R}^d}V_n^{\mathrm{out}}|\zeta_nu|^2\,dx
	+
	\int_{A_n}V_n|u|^2\,dx.
	\]
	Combining this identity with \eqref{eq:splitting-annulus-bound}, we have
	\[
	q_{V_n}[u]
	\geq
	(1-\varepsilon_n)\|(-\Delta)^{s/2}u\|_2^2
	-\int_{\mathbb{R}^d}V_n^{\mathrm{in}}|\chi_nu|^2\,dx
	-\int_{\mathbb{R}^d}V_n^{\mathrm{out}}|\zeta_nu|^2\,dx
	-\varepsilon_n\|u\|_2^2.
	\]
	Accordingly, this yields
	\begin{equation}\label{eq:splitting-form}
		\begin{aligned}
			q_{V_n}[u]\geq{}&
			(1-\varepsilon_n)
			\sum_{h\in\{\chi_n,\zeta_n,\eta_n\}}
			\|(-\Delta)^{s/2}(hu)\|_2^2
			-\int_{\mathbb{R}^d}V_n^{\mathrm{in}}|\chi_nu|^2\,dx\\
			&-\int_{\mathbb{R}^d}V_n^{\mathrm{out}}|\zeta_nu|^2\,dx
			-b_n\|u\|_2^2,
		\end{aligned}
	\end{equation}
	where $b_n:=C(R_n^{-2s}+\varepsilon_n)\to0$ since $R_n\to\infty$ and $\varepsilon_n\to0$.
	To express the preceding lower bound in operator form, we define, for $c\in[1/2,1]$,
	\[
	D_{n,c}:=
	\bigl(c(-\Delta)^s-V_n^{\mathrm{in}}\bigr)
	\oplus
	\bigl(c(-\Delta)^s-V_n^{\mathrm{out}}\bigr)
	\oplus
	c(-\Delta)^s.
	\]
	Since $\chi_n^2+\zeta_n^2+\eta_n^2=1$, the map
	\[
	J_nu:=(\chi_nu,\zeta_nu,\eta_nu)
	\]
	is an $L^2$-isometry from $H^s(\mathbb{R}^d)$ into the form domain of $D_{n,c}$. With this notation, \eqref{eq:splitting-form} can be written as
	\[
	q_{V_n}[u]
	\geq
	q_{D_{n,1-\varepsilon_n}}[J_nu]
	-b_n\|u\|_2^2.
	\]
	By applying the min--max principle, we obtain
	\begin{equation}\label{eq:split-preliminary}
		\lambda_j(H_{V_n})
		\geq
		\lambda_j(D_{n,1-\varepsilon_n})-b_n,
		\qquad 1\leq j\leq N.
	\end{equation}
	
	Since $\{V_n\}$ is bounded in $L^p(\mathbb{R}^d)$, fix $M>0$ such that 
	$$\sup_n\|V_n\|_p\leq M.$$ 
	The bounds
	$\|V_n^{\mathrm{in}}\|_p,\|V_n^{\mathrm{out}}\|_p\leq M$
	allow us to apply \eqref{eq:kinetic-continuity} of
	Lemma~\ref{lem:spectral-stability} to the first two components of
	$D_{n,c}$, while the third component $c(-\Delta)^s$ has no negative
	spectrum. The min--max principle then yields
	\begin{equation}\label{eq:direct-sum-kinetic}
		0\leq
		\lambda_j(D_{n,1})-\lambda_j(D_{n,c})
		\leq
		C_M(1-c),
		\qquad 1\leq j\leq N.
	\end{equation}
	Furthermore, define
	\[
	B_n:=H_{V_n^{\mathrm{in}}}\oplus H_{V_n^{\mathrm{out}}},
	\qquad
	\mu_{j,n}:=\lambda_j(B_n),
	\qquad
	r_n:=b_n+C_M\varepsilon_n.
	\]
	Then $\lambda_j(D_{n,1})=\mu_{j,n}$ and, since $b_n\to0$ and $\varepsilon_n\to0$, we have $r_n\to0$.
	Taking
	$c=1-\varepsilon_n$ in \eqref{eq:direct-sum-kinetic}, we deduce from
	\eqref{eq:split-preliminary} that
	\begin{equation}\label{eq:split-levels}
		\lambda_j(H_{V_n})
		\geq
		\lambda_j(D_{n,1-\varepsilon_n})-b_n
		\ge \lambda_j(D_{n,1})-C_M\varepsilon_n-b_n
		=
		\mu_{j,n}-r_n,
		\qquad 1\leq j\leq N.
	\end{equation}
	This completes the lower spectral estimate.
	
	\textbf{Step 2. Truncation and interaction estimates.}
	We first truncate eigenfunctions corresponding to spectral levels in
	$(-\infty,-E_0]$, where $E_0>0$ is fixed, so that the interior and
	exterior truncated functions have separated supports.
	Take
	$\theta_n^{\mathrm{in}},\theta_n^{\mathrm{out}}
	\in C^\infty(\mathbb{R}^d;[0,1])$ such that
	\[
	\left\{
	\begin{aligned}
		\theta_n^{\mathrm{in}}&=0 &&\text{on }B_{R_n},\\
		\theta_n^{\mathrm{in}}&=1 &&\text{on }B_{5R_n/4}^c,
	\end{aligned}
	\right.
	\qquad
	\left\{
	\begin{aligned}
		\theta_n^{\mathrm{out}}&=1 &&\text{on }B_{7R_n/4},\\
		\theta_n^{\mathrm{out}}&=0 &&\text{on }B_{2R_n}^c,
	\end{aligned}
	\right.
	\qquad
	\|\nabla\theta_n^{\mathrm{in}}\|_\infty
	+\|\nabla\theta_n^{\mathrm{out}}\|_\infty
	\leq \frac{C}{R_n}.
	\]
	Let $W_n$ denote either $V_n^{\mathrm{in}}$ or $V_n^{\mathrm{out}}$, with
	$\theta_n$ the corresponding cutoff, namely, 
	$$(W_n,\theta_n)\in
	\left\{
	\bigl(V_n^{\mathrm{in}},\theta_n^{\mathrm{in}}\bigr),
	\bigl(V_n^{\mathrm{out}},\theta_n^{\mathrm{out}}\bigr)
	\right\}.$$
	Obviously, $W_n\theta_n=0$ almost everywhere. 
	Suppose that $u_n$ is an $L^2$-normalized eigenfunction of $H_{W_n}$ with eigenvalue $\lambda_n\leq-E_0$, that is,
	\begin{equation}\label{eq:eigenvalue-equation-Wn}
		H_{W_n}u_n=(-\Delta)^su_n-W_nu_n=\lambda_nu_n.
	\end{equation}
	We claim that $\theta_n^2u_n\in H^s(\mathbb{R}^d)$. First,
	$\theta_n^2u_n\in L^2(\mathbb{R}^d)$ since $0\leq\theta_n\leq1$. Moreover, a direct computation gives
	\[
	\begin{aligned}
		|\theta_n(x)^2u_n(x)-\theta_n(y)^2u_n(y)|^2
		&=
		\left|
		\theta_n(x)^2\bigl(u_n(x)-u_n(y)\bigr)
		+
		u_n(y)\bigl(\theta_n(x)^2-\theta_n(y)^2\bigr)
		\right|^2\\
		&\quad\leq
		2\theta_n(x)^4|u_n(x)-u_n(y)|^2
		+
		2|u_n(y)|^2
		|\theta_n(x)^2-\theta_n(y)^2|^2\\
		&\quad\leq
		2|u_n(x)-u_n(y)|^2
		+
		8|u_n(y)|^2|\theta_n(x)-\theta_n(y)|^2.
	\end{aligned}
	\]
	Noting that
	$|\theta_n(x)-\theta_n(y)|^2
	\leq
	C\min\left\{1,\frac{|x-y|^2}{R_n^2}\right\},$
	and arguing as in the proof of
	Lemma~\ref{lem:fractional-localization}-\ref{item:discrete-localization} (see \eqref{eq:IMS-kernel-identity}--\eqref{eq:IMS-kernel-scaling}),
	we obtain
	\[
	\iint_{\mathbb{R}^d\times\mathbb{R}^d}
	\frac{|\theta_n(x)^2u_n(x)-\theta_n(y)^2u_n(y)|^2}
	{|x-y|^{d+2s}}\,dx\,dy
	<\infty.
	\]
	Hence $\theta_n^2u_n\in H^s(\mathbb{R}^d)$.
	Testing the weak form of \eqref{eq:eigenvalue-equation-Wn} with
	$\theta_n^2u_n$ and using $W_n\theta_n=0$, we obtain
	\begin{equation}\label{eq:cutoff-eigenvalue}
		\operatorname{Re}\int_{\mathbb{R}^d}
		(-\Delta)^{s/2}u_n\,
		\overline{(-\Delta)^{s/2}(\theta_n^2u_n)}\,dx
		=
		\lambda_n\|\theta_nu_n\|_2^2.
	\end{equation}
	Moreover, by the kernel representation \eqref{eq:kernel-form}, we note that
	\[
	\begin{aligned}
		\|(-\Delta)^{s/2}(\theta_nu_n)\|_2^2
		-\operatorname{Re}\int_{\mathbb{R}^d}
		(-\Delta)^{s/2}u_n\,
		\overline{(-\Delta)^{s/2}(\theta_n^2u_n)}\,dx=
		\frac{c_{d,s}}2
		\iint_{\mathbb{R}^d\times\mathbb{R}^d}
		\frac{\mathcal I_n(x,y)}
		{|x-y|^{d+2s}}\,dx\,dy,
	\end{aligned}
	\]
	where
	\[
	\begin{aligned}
		\mathcal I_n(x,y)
		&=
		|\theta_n(x)u_n(x)-\theta_n(y)u_n(y)|^2-
		\operatorname{Re}\!\left[
		(u_n(x)-u_n(y))
		\overline{\theta_n(x)^2u_n(x)-\theta_n(y)^2u_n(y)}
		\right]\\
		&=
		\theta_n(x)^2|u_n(x)|^2
		+\theta_n(y)^2|u_n(y)|^2
		-2\theta_n(x)\theta_n(y)
		\operatorname{Re}(u_n(x)\overline{u_n(y)})\\
		&\quad-
		\theta_n(x)^2|u_n(x)|^2
		-\theta_n(y)^2|u_n(y)|^2
		+\bigl(\theta_n(x)^2+\theta_n(y)^2\bigr)
		\operatorname{Re}(u_n(x)\overline{u_n(y)})\\
		&=
		|\theta_n(x)-\theta_n(y)|^2
		\operatorname{Re}(u_n(x)\overline{u_n(y)}).
	\end{aligned}
	\]
	It thus follows from \eqref{eq:cutoff-eigenvalue} that
	\[
	\|(-\Delta)^{s/2}(\theta_nu_n)\|_2^2
	-\lambda_n\|\theta_nu_n\|_2^2=
	\frac{c_{d,s}}2
	\iint_{\mathbb{R}^d\times\mathbb{R}^d}
	\frac{|\theta_n(x)-\theta_n(y)|^2
		\operatorname{Re}(u_n(x)\overline{u_n(y)})}
	{|x-y|^{d+2s}}\,dx\,dy.
	\]
	The same argument as in \eqref{eq:IMS-kernel-identity}--\eqref{eq:IMS-kernel-scaling} gives
	\[
	\|(-\Delta)^{s/2}(\theta_nu_n)\|_2^2
	-\lambda_n\|\theta_nu_n\|_2^2
	\leq
	CR_n^{-2s}\|u_n\|_2^2
	=
	CR_n^{-2s},
	\]
	which, together with the fact that $\lambda_n\leq-E_0$, indicates that
	\begin{equation}\label{eq:splitting-tail-bound}
		\|(-\Delta)^{s/2}(\theta_nu_n)\|_2^2
		+E_0\|\theta_nu_n\|_2^2\leq CR_n^{-2s}.
	\end{equation}
	Set $$h_n:=(1-\theta_n)u_n.$$ 
	By \eqref{eq:splitting-tail-bound}, one has
	\begin{equation}\label{eq:splitting-truncation}
		\|h_n-u_n\|_{H^s}
		=
		\|\theta_nu_n\|_{H^s}
		\leq C_{E_0}R_n^{-s}.
	\end{equation}
	Moreover, \eqref{eq:uniform-kinetic} with $a=1$ implies
	$\|u_n\|_{H^s}\leq C$, and hence 
	\begin{equation}\label{eq:truncation-uniform-bound}
		\|h_n\|_{H^s}\leq C_{E_0}.
	\end{equation}
	By construction, for the interior component, $h_n$ is supported in $\overline{B}_{5R_n/4}$ and coincides with $u_n$ on $B_{R_n}$; for the exterior component, $h_n$ vanishes on $B_{7R_n/4}$ and coincides with $u_n$ on $B_{2R_n}^c$. 
	
	We next estimate the interaction between the two truncated components.
	Let $f,g\in H^s(\mathbb{R}^d)$ be supported in $\overline{B}_{5R_n/4}$ and $B_{7R_n/4}^c$, respectively. Since their supports are separated by a distance of at least $R_n/2$, it follows that
	\[
	\int_{\mathbb{R}^d}
	(-\Delta)^{s/2}f\,
	\overline{(-\Delta)^{s/2}g}\,dx
	=
	-c_{d,s}
	\int_{B_{5R_n/4}}
	\int_{B_{7R_n/4}^c}
	\frac{f(x)\overline{g(y)}}{|x-y|^{d+2s}}\,dy\,dx.
	\]
	Letting $w:=x-y$, we have $|w|\geq R_n/2$ whenever $f(x)g(y)\neq0$.
	Hence, by the Cauchy--Schwarz inequality and Tonelli's theorem,
	\begin{equation}\label{eq:splitting-cross-term}
		\left|
		\int_{\mathbb{R}^d}
		(-\Delta)^{s/2}f\,
		\overline{(-\Delta)^{s/2}g}\,dx
		\right|
		\leq
		c_{d,s}\|f\|_2\|g\|_2
		\int_{|w|\geq R_n/2}|w|^{-d-2s}\,dw
		\leq
		CR_n^{-2s}\|f\|_2\|g\|_2.
	\end{equation}
	
	\textbf{Step 3. Gram matrix and energy matrix estimates.}
	Set
	\[
	K_n:=\#\{1\leq j\leq N:\mu_{j,n}\leq-E_0\}.
	\]
	That is, $K_n$ is the number of the first $N$ min--max levels of $B_n$
	lying in $(-\infty,-E_0]$. If $K_n=0$, there is nothing to prove in
	this step, so we assume that $K_n\geq1$. 
	By the direct-sum structure  $B_n=H_{V_n^{\mathrm{in}}}\oplus H_{V_n^{\mathrm{out}}}$, the corresponding eigenvectors can be chosen real and orthonormal, each of the form $(u,0)$ or $(0,v)$. Denote by $u_{1,n},\ldots,u_{K_n,n}$ the $L^2$-normalized nonzero components of these eigenvectors. 
	For $r=1,\ldots,K_n$, define
	\[
	\theta_{r,n}:=
	\begin{cases}
		\theta_n^{\mathrm{in}}, & \text{if }u_{r,n}\text{ corresponds to }
		H_{V_n^{\mathrm{in}}},\\
		\theta_n^{\mathrm{out}}, & \text{if }u_{r,n}\text{ corresponds to }
		H_{V_n^{\mathrm{out}}}.
	\end{cases}
	\]
	Let $h_{r,n}:=(1-\theta_{r,n})u_{r,n}$. 
	Applying \eqref{eq:splitting-truncation} to each $u_{r,n}$, we obtain
	\begin{equation}\label{eq:splitting-truncated-functions}
		\|h_{r,n}-u_{r,n}\|_{H^s}
		\leq C_{E_0}R_n^{-s},
		\qquad
		1\leq r\leq K_n.
	\end{equation}
	Moreover, \eqref{eq:uniform-kinetic} with $a=1$ and
	\eqref{eq:truncation-uniform-bound} imply
	\begin{equation}\label{eq:splitting-uniform-Hs-bound}
		\|u_{r,n}\|_{H^s}+\|h_{r,n}\|_{H^s}
		\leq C_{E_0},
		\qquad
		1\leq r\leq K_n.
	\end{equation}
	For $h_{r,n}$ and $h_{t,n}$ arising from the same component,
	we deduce from \eqref{eq:splitting-truncated-functions} and \eqref{eq:splitting-uniform-Hs-bound} that
	\[
	\left|\langle h_{r,n},h_{t,n}\rangle-\delta_{rt}\right|
	\leq
	\|h_{r,n}-u_{r,n}\|_2\|h_{t,n}\|_2
	+\|u_{r,n}\|_2\|h_{t,n}-u_{t,n}\|_2
	\leq
	C_{E_0}R_n^{-s}.
	\]
	On the other hand, if they arise from different components, their supports are disjoint, and hence $\langle h_{r,n},h_{t,n}\rangle=0$ and $\delta_{rt}=0$. Therefore, the Gram matrix
	\begin{equation*}
		G_n:=
		\bigl(\langle h_{r,n},h_{t,n}\rangle\bigr)_{r,t=1}^{K_n}
	\end{equation*}
	satisfies, since $K_n\leq N$,
	\begin{equation}\label{eq:splitting-Gram}
		\|G_n-I\|_{\mathrm{op}}
		\leq
		K_n\max_{1\leq r,t\leq K_n}
		\left|
		\langle h_{r,n},h_{t,n}\rangle-\delta_{rt}
		\right| 
		\leq C_{E_0}R_n^{-s}\to0,
	\end{equation}
	where $\|\cdot\|_{\mathrm{op}}$ denotes the operator norm on
	$\mathbb{R}^{K_n}$. This proves the Gram matrix estimate.
	
	We next turn to the energy matrix estimates. Define the energy matrix $\mathcal A_n$ and the diagonal matrix
	$\Lambda_n$ by
	\begin{equation}\label{eq:splitting-energy-definition}
		(\mathcal A_n)_{rt}
		:=
		\langle(-\Delta)^{s/2}h_{r,n},(-\Delta)^{s/2}h_{t,n}\rangle
		-\int_{\mathbb R^d}V_nh_{r,n}h_{t,n}\,dx,
		\quad
		\Lambda_n:=\operatorname{diag}(\mu_{1,n},\ldots,\mu_{K_n,n}).
	\end{equation}
	We first consider the case where $u_{r,n}$ and $u_{t,n}$ arise from the same component. By \eqref{eq:splitting-truncated-functions} and \eqref{eq:splitting-uniform-Hs-bound}, we have 
	\begin{equation}\label{eq:splitting-kinetic-error}
		\begin{aligned}
			&\bigl|
			\langle(-\Delta)^{s/2}h_{r,n},(-\Delta)^{s/2}h_{t,n}\rangle
			-\langle(-\Delta)^{s/2}u_{r,n},(-\Delta)^{s/2}u_{t,n}\rangle
			\bigr| \\
			&\leq
			\|(-\Delta)^{s/2}(h_{r,n}-u_{r,n})\|_2
			\|(-\Delta)^{s/2}h_{t,n}\|_2 +
			\|(-\Delta)^{s/2}u_{r,n}\|_2
			\|(-\Delta)^{s/2}(h_{t,n}-u_{t,n})\|_2\\
			&\leq C_{E_0}R_n^{-s}.
		\end{aligned}
	\end{equation}
	If they both arise from $H_{V_n^{\mathrm{in}}}$, then
	$\theta_{r,n}=\theta_{t,n}=\theta_n^{\mathrm{in}}$, and hence
	$V_n^{\mathrm{in}}\theta_{r,n}
	=V_n^{\mathrm{in}}\theta_{t,n}=0$ a.e. in $\mathbb R^d$. Therefore,
	\begin{equation*}
		\int_{\mathbb R^d}V_n^{\mathrm{in}}h_{r,n}h_{t,n}\,dx
		=
		\int_{\mathbb R^d}
		V_n^{\mathrm{in}}(1-\theta_{r,n})u_{r,n}
		(1-\theta_{t,n})u_{t,n}\,dx
		=
		\int_{\mathbb R^d}V_n^{\mathrm{in}}u_{r,n}u_{t,n}\,dx.
	\end{equation*}
	Moreover,
	$\operatorname{supp}h_{r,n},\operatorname{supp}h_{t,n}
	\subset\overline{B_{5R_n/4}}$ and
	$\operatorname{supp}V_n^{\mathrm{out}}\subset B_{2R_n}^c$, so that
	\begin{equation*}
		\int_{\mathbb R^d}V_n^{\mathrm{out}}h_{r,n}h_{t,n}\,dx=0.
	\end{equation*}
	If $u_{r,n}$ and $u_{t,n}$ both arise from
	$H_{V_n^{\mathrm{out}}}$, then similarly
	$\int_{\mathbb R^d}V_n^{\mathrm{out}}h_{r,n}h_{t,n}\,dx
	=
	\int_{\mathbb R^d}V_n^{\mathrm{out}}u_{r,n}u_{t,n}\,dx,$
	while
	$\int_{\mathbb R^d}V_n^{\mathrm{in}}h_{r,n}h_{t,n}\,dx=0.$
	In either case, by H\"older's inequality, \eqref{eq:GN},
	\eqref{eq:splitting-uniform-Hs-bound}, and
	\eqref{eq:splitting-annulus-notation}, there holds
	\begin{equation*}
		\left|\int_{A_n}V_nh_{r,n}h_{t,n}\,dx\right|\leq
		\|V_n\mathbf 1_{A_n}\|_p
		\|h_{r,n}\|_{2p'}\|h_{t,n}\|_{2p'}
		\leq C_{E_0}\varepsilon_n.
	\end{equation*}
	Furthermore, since $u_{r,n}$ and $u_{t,n}$ are orthonormal eigenfunctions of the same operator, when they both arise from $H_{V_n^{\mathrm{in}}}$, we have
	\begin{equation}\label{eq:splitting-component-eigenvalue}
		\langle(-\Delta)^{s/2}u_{r,n},(-\Delta)^{s/2}u_{t,n}\rangle
		-\int_{\mathbb R^d}V_n^{\mathrm{in}}u_{r,n}u_{t,n}\,dx
		=
		\mu_{r,n}\delta_{rt},
	\end{equation}
	and the analogous identity holds for $H_{V_n^{\mathrm{out}}}$. When $u_{r,n}$ and $u_{t,n}$ both arise from
	$H_{V_n^{\mathrm{in}}}$, it follows from
	\eqref{eq:splitting-kinetic-error}--\eqref{eq:splitting-component-eigenvalue} that
	\begin{equation*}
		\begin{aligned}
			\left|(\mathcal A_n)_{rt}-\mu_{r,n}\delta_{rt}\right|
			&=
			\Bigg|
			\langle(-\Delta)^{s/2}h_{r,n},(-\Delta)^{s/2}h_{t,n}\rangle
			-\int_{\mathbb R^d}V_nh_{r,n}h_{t,n}\,dx \\
			&\qquad
			-\langle(-\Delta)^{s/2}u_{r,n},(-\Delta)^{s/2}u_{t,n}\rangle
			+\int_{\mathbb R^d}V_n^{\mathrm{in}}u_{r,n}u_{t,n}\,dx
			\Bigg| \\
			&\leq
			\left|
			\langle(-\Delta)^{s/2}h_{r,n},(-\Delta)^{s/2}h_{t,n}\rangle
			-\langle(-\Delta)^{s/2}u_{r,n},(-\Delta)^{s/2}u_{t,n}\rangle
			\right| 
			+\left|\int_{A_n}V_nh_{r,n}h_{t,n}\,dx\right| \\
			&\leq C_{E_0}\bigl(R_n^{-s}+\varepsilon_n\bigr).
		\end{aligned}
	\end{equation*}
	If they both arise from $H_{V_n^{\mathrm{out}}}$, the same argument,
	with $V_n^{\mathrm{out}}$ in place of $V_n^{\mathrm{in}}$, gives the
	same estimate. Hence, whenever $u_{r,n}$ and $u_{t,n}$ arise from the
	same component,
	\begin{equation}\label{eq:splitting-energy-entry}
		\left|(\mathcal A_n)_{rt}-\mu_{r,n}\delta_{rt}\right|
		\leq
		C_{E_0}\bigl(R_n^{-s}+\varepsilon_n\bigr).
	\end{equation}

	It remains to consider the case in which $u_{r,n}$ and $u_{t,n}$
	arise from different components.
	Up to interchanging $r$ and $t$, assume that $u_{r,n}$ comes from the interior component and $u_{t,n}$ from the exterior component. Then
	$\operatorname{supp}h_{r,n}\subset \overline{B_{5R_n/4}},
	\operatorname{supp}h_{t,n}\subset B_{7R_n/4}^c.$
	Hence their supports are disjoint, and
	\begin{equation}\label{eq:splitting-cross-potential}
		\int_{\mathbb R^d}V_nh_{r,n}h_{t,n}\,dx=0.
	\end{equation}
	Moreover, by \eqref{eq:splitting-cross-term} and \eqref{eq:splitting-uniform-Hs-bound}, we get
	\begin{equation}\label{eq:splitting-cross-kinetic}
		\left|
		\langle(-\Delta)^{s/2}h_{r,n},
		(-\Delta)^{s/2}h_{t,n}\rangle
		\right|
		\leq C_{E_0}R_n^{-2s}.
	\end{equation}
	For $u_{r,n}$ and $u_{t,n}$ arising from different components, we have $r\neq t$ and hence $\delta_{rt}=0$. Therefore, it follows from
	\eqref{eq:splitting-cross-potential} and
	\eqref{eq:splitting-cross-kinetic} that
	\begin{equation}\label{eq:splitting-cross-energy-entry}
		\left|(\mathcal A_n)_{rt}-\mu_{r,n}\delta_{rt}\right|
		=
		\left|
		\langle(-\Delta)^{s/2}h_{r,n},
		(-\Delta)^{s/2}h_{t,n}\rangle
		\right| 
		\leq C R_n^{-2s}
		\leq C R_n^{-s},
	\end{equation}
	where the last inequality follows from $R_n\geq1$ and $s>0$.
	
	Consequently, by \eqref{eq:splitting-energy-entry} and
	\eqref{eq:splitting-cross-energy-entry}, we obtain
	\[
	\left|(\mathcal A_n)_{rt}-\mu_{r,n}\delta_{rt}\right|
	\leq C_{E_0}\bigl(R_n^{-s}+\varepsilon_n\bigr),
	\qquad 1\leq r,t\leq K_n.
	\]
	Since $K_n\leq N$, it follows that
	\begin{equation}\label{eq:splitting-energy-matrix}
		\|\mathcal A_n-\Lambda_n\|_{\mathrm{op}}
		\leq C_{E_0}\bigl(R_n^{-s}+\varepsilon_n\bigr)\to0.
	\end{equation}

	\textbf{Step 4. Orthogonalization and upper spectral estimate.}
	As in Step~3, it suffices to consider the case $K_n\geq1$. 
	For all sufficiently large $n$, \eqref{eq:splitting-Gram} gives
	$\|G_n-I\|_{\mathrm{op}}\leq1/2$. Since $G_n$ is symmetric, its spectrum satisfies
	\[
	\sigma(G_n)
	\subset
	[1-\|G_n-I\|_{\mathrm{op}},\,1+\|G_n-I\|_{\mathrm{op}}]
	\subset[1/2,3/2].
	\]
	Hence, $G_n$ is positive definite and $G_n^{-1/2}$ is well defined.
	For $1\leq r\leq K_n$, set
	\begin{equation}\label{eq:splitting-orthogonalized-basis}
		\widetilde h_{r,n}:=
		\sum_{\ell=1}^{K_n}h_{\ell,n}(G_n^{-1/2})_{\ell r}.
	\end{equation}
	Since $G_n^{-1/2}$ is symmetric, for $1\leq r,t\leq K_n$ we have
	\[
	\begin{aligned}
		\langle \widetilde h_{r,n},\widetilde h_{t,n}\rangle
		&=
		\left\langle
		\sum_{a=1}^{K_n}h_{a,n}(G_n^{-1/2})_{ar},
		\sum_{b=1}^{K_n}h_{b,n}(G_n^{-1/2})_{bt}
		\right\rangle\\
		&=
		\sum_{a,b=1}^{K_n}
		(G_n^{-1/2})_{ar}
		\langle h_{a,n},h_{b,n}\rangle
		(G_n^{-1/2})_{bt} \\
		&=
		\bigl(G_n^{-1/2}G_nG_n^{-1/2}\bigr)_{rt}
		=\delta_{rt}.
	\end{aligned}
	\]
	Thus, $\widetilde h_{1,n},\ldots,\widetilde h_{K_n,n}$ are orthonormal in $L^2(\mathbb R^d)$.
	
	Let $\widetilde{\mathcal A}_n$ denote the energy matrix of the orthonormal family
	$\widetilde h_{1,n},\ldots,\widetilde h_{K_n,n}$, that is,
	\[
	(\widetilde{\mathcal A}_n)_{rt}
	:=
	q_{V_n}[\widetilde h_{r,n},\widetilde h_{t,n}]
	=
	\langle(-\Delta)^{s/2}\widetilde h_{r,n},
	(-\Delta)^{s/2}\widetilde h_{t,n}\rangle
	-\int_{\mathbb R^d}V_n\widetilde h_{r,n}\widetilde h_{t,n}\,dx.
	\]	
	By \eqref{eq:splitting-energy-definition} and
	\eqref{eq:splitting-orthogonalized-basis}, we get
	\[
	\widetilde{\mathcal A}_n
	=
	G_n^{-1/2}\mathcal A_nG_n^{-1/2}.
	\]
	For $n$ large enough, since $G_n$ is symmetric positive definite and
	$\sigma(G_n)\subset[1/2,3/2]$, the spectral theorem yields
	\begin{equation}\label{eq:splitting-Gram-inverse-bound}
		\|G_n^{-1/2}\|_{\mathrm{op}}
		=
		\max_{\lambda\in\sigma(G_n)}\lambda^{-1/2}
		\leq\sqrt2.
	\end{equation}
	Furthermore, from \eqref{eq:splitting-Gram} and the Lipschitz continuity of $x\mapsto x^{-1/2}$ on $[1/2,3/2]$, we deduce that
	\begin{equation}\label{eq:splitting-Gram-inverse}
		\|G_n^{-1/2}-I\|_{\mathrm{op}}
		=
		\max_{\lambda\in\sigma(G_n)}
		|\lambda^{-1/2}-1| \\
		\leq
		C\max_{\lambda\in\sigma(G_n)}|\lambda-1| 
		=
		C\|G_n-I\|_{\mathrm{op}}
		\to0.		
	\end{equation}
	Recall from \eqref{eq:splitting-energy-definition} that
	$\Lambda_n=\operatorname{diag}(\mu_{1,n},\ldots,\mu_{K_n,n})$.
	By \eqref{eq:ground-bound} and the uniform $L^p$ bounds on
	$V_n^{\mathrm{in}}$ and $V_n^{\mathrm{out}}$, we obtain
	\begin{equation}\label{eq:splitting-Lambda-bound}
		\|\Lambda_n\|_{\mathrm{op}}
		=
		\max_{1\leq j\leq K_n}|\mu_{j,n}|
		=
		|\mu_{1,n}|
		\leq C.
	\end{equation}
	Observe that
	\[
	\widetilde{\mathcal A}_n-\Lambda_n
	=G_n^{-1/2}(\mathcal A_n-\Lambda_n)G_n^{-1/2}
	+(G_n^{-1/2}-I)\Lambda_nG_n^{-1/2}
	+\Lambda_n(G_n^{-1/2}-I).
	\]
	Therefore, by \eqref{eq:splitting-energy-matrix}, \eqref{eq:splitting-Gram-inverse-bound},
	\eqref{eq:splitting-Gram-inverse}, and
	\eqref{eq:splitting-Lambda-bound}, one has 
	\[
	\begin{aligned}
		\|\widetilde{\mathcal A}_n-\Lambda_n\|_{\mathrm{op}}
		&\leq
		\|G_n^{-1/2}\|_{\mathrm{op}}^2
		\|\mathcal A_n-\Lambda_n\|_{\mathrm{op}} 
		+
		\|G_n^{-1/2}-I\|_{\mathrm{op}}
		\|\Lambda_n\|_{\mathrm{op}}
		\|G_n^{-1/2}\|_{\mathrm{op}} 
		+
		\|\Lambda_n\|_{\mathrm{op}}
		\|G_n^{-1/2}-I\|_{\mathrm{op}} \\
		&\leq
		C\bigl(
		\|\mathcal A_n-\Lambda_n\|_{\mathrm{op}}
		+\|G_n-I\|_{\mathrm{op}}
		\bigr)
		\to0
		\qquad\text{as } n\to\infty.
	\end{aligned}
	\]
	Denote by $\nu_{1,n}\leq\cdots\leq\nu_{K_n,n}$ the eigenvalues of
	$\widetilde{\mathcal A}_n$. Noting that
	$\mu_{1,n},\ldots,\mu_{K_n,n}$ are the eigenvalues of
	$\Lambda_n$, the finite-dimensional min--max principle gives
	\begin{equation}\label{eq:splitting-finite-dimensional}
		|\nu_{j,n}-\mu_{j,n}|
		\leq
		\|\widetilde{\mathcal A}_n-\Lambda_n\|_{\mathrm{op}}
		\to0,
		\qquad 1\leq j\leq K_n.
	\end{equation}
	Since $\mu_{j,n}\leq-E_0$ for $1\leq j\leq K_n$ and
	$\|\widetilde{\mathcal A}_n-\Lambda_n\|_{\mathrm{op}}\leq E_0/2$
	for $n$ large enough, it therefore follows that
	\begin{equation}\label{eq:splitting-nu-negative}
		\nu_{j,n}
		\leq
		\mu_{j,n}
		+\|\widetilde{\mathcal A}_n-\Lambda_n\|_{\mathrm{op}}
		\leq
		-\frac{E_0}{2}<0,
		\qquad 1\leq j\leq K_n.
	\end{equation}
	
	Finally, set
	\[
	\widetilde E_n:=
	\operatorname{span}\{\widetilde h_{1,n},\ldots,\widetilde h_{K_n,n}\}.
	\]
	For each $1\leq j\leq K_n$, let $\widetilde E_{j,n}\subset\widetilde E_n$
	be the $j$-dimensional subspace corresponding, under the basis
	$\{\widetilde h_{r,n}\}_{r=1}^{K_n}$, to the span of orthonormal eigenvectors of $\widetilde{\mathcal A}_n$ associated with
	$\nu_{1,n},\ldots,\nu_{j,n}$. Then the restriction of $q_{V_n}$ to
	$\widetilde E_{j,n}$ is represented in this eigenbasis by
	$\operatorname{diag}(\nu_{1,n},\ldots,\nu_{j,n})$.	
	Writing $u\in\widetilde E_{j,n}$ in the corresponding orthonormal
	eigenbasis with coordinates $b_1,\ldots,b_j$, we have
	$\sum_{\ell=1}^j|b_\ell|^2=\|u\|_2^2$
	and
	\[
	\sup_{\substack{u\in\widetilde E_{j,n}\\ \|u\|_2=1}}
	q_{V_n}[u]
	=
	\sup_{\sum_{\ell=1}^j|b_\ell|^2=1}
	\sum_{\ell=1}^j\nu_{\ell,n}|b_\ell|^2
	=
	\nu_{j,n}.
	\]
	For sufficiently large $n$, applying the min--max principle for $\lambda_j(H_{V_n})$ with
	$\widetilde E_{j,n}$ as a trial space and using
	\eqref{eq:splitting-nu-negative}, we obtain
	\[
	\lambda_j(H_{V_n})
	\leq
	\sup_{\substack{u\in\widetilde E_{j,n}\\ \|u\|_2=1}}
	q_{V_n}[u]
	=
	\nu_{j,n}<0,
	\qquad 1\leq j\leq K_n.
	\]
	Together with \eqref{eq:splitting-finite-dimensional}, this yields
	\begin{equation}\label{eq:split-levels-upper}
		\lambda_j(H_{V_n})
		\leq \nu_{j,n}
		\leq \mu_{j,n}
		+\|\widetilde{\mathcal A}_n-\Lambda_n\|_{\mathrm{op}},
		\qquad 1\leq j\leq K_n.
	\end{equation}
	This completes the upper spectral estimate.
	
	\textbf{Step 5. Conclusion.}
	Fix $E_0>0$. For $1\leq j\leq K_n$, combining \eqref{eq:split-levels} with
	\eqref{eq:split-levels-upper}, we have
	\begin{equation}\label{eq:splitting-level-convergence}
		|\lambda_j(H_{V_n})-\mu_{j,n}|
		\leq
		\max\left\{
		r_n,\,
		\|\widetilde{\mathcal A}_n-\Lambda_n\|_{\mathrm{op}}
		\right\}
		\to0
		\qquad\text{as }n\to\infty.
	\end{equation}
	By \eqref{eq:ground-bound} and the uniform $L^p$ bounds, there exists $C>0$, independent of $n$, such that
	\begin{equation*}\label{eq:splitting-uniform-spectral-bound}
		-C\leq \lambda_j(H_{V_n})\leq0,
		\qquad
		-C\leq \mu_{j,n}\leq0,
		\qquad 1\leq j\leq N.
	\end{equation*}	
	Consequently, by the uniform continuity of $t\mapsto |t|^\gamma$ on $[-C,0]$ and using $K_n\leq N$, it follows from \eqref{eq:splitting-level-convergence} that
	\[
	\sum_{j=1}^{K_n}
	\left|
	|\lambda_j(H_{V_n})|^\gamma-|\mu_{j,n}|^\gamma
	\right|
	\to0
	\qquad\text{as }n\to\infty.
	\]
	On the other hand, for $K_n<j\leq N$, the definition of $K_n$ gives
	$-E_0<\mu_{j,n}\leq0$, while \eqref{eq:split-levels} yields $-E_0-r_n<\lambda_j(H_{V_n})\leq0.$
	Thus,
	\[
	\sum_{j=K_n+1}^N
	\left|
	|\lambda_j(H_{V_n})|^\gamma-|\mu_{j,n}|^\gamma
	\right|
	\leq
	N\bigl((E_0+r_n)^\gamma+E_0^\gamma\bigr).
	\]
	Combining the two estimates and letting $n\to\infty$, one can see that
	\[
	\limsup_{n\to\infty}
	\left|
	\sum_{j=1}^N|\lambda_j(H_{V_n})|^\gamma
	-\sum_{j=1}^N|\mu_{j,n}|^\gamma
	\right|
	\leq 2NE_0^\gamma.
	\]
	Furthermore, letting $E_0\to0$ and using $\gamma>0$, we conclude that
	\begin{equation}\label{eq:splitting-direct-sum}
		\sum_{j=1}^N|\lambda_j(H_{V_n})|^\gamma
		=
		\sum_{j=1}^N|\mu_{j,n}|^\gamma+o_n(1).
	\end{equation}
	Recall that
	$B_n=H_{V_n^{\mathrm{in}}}\oplus H_{V_n^{\mathrm{out}}}$, whose first
	$N$ min--max levels are obtained from those of the two components.
	Hence, there exists $k_n\in\{0,\ldots,N\}$ such that
	\begin{equation}\label{eq:splitting-direct-sum-levels}
		\sum_{j=1}^N|\mu_{j,n}|^\gamma
		=
		\sum_{j=1}^{k_n}|\lambda_j(H_{V_n^{\mathrm{in}}})|^\gamma
		+
		\sum_{j=1}^{N-k_n}|\lambda_j(H_{V_n^{\mathrm{out}}})|^\gamma.
	\end{equation}
	Since $k_n\in\{0,\ldots,N\}$, after passing to a subsequence we may assume that $k_n=k$ for some $k\in\{0,\ldots,N\}$. It then follows from \eqref{eq:splitting-direct-sum-levels} and
	\eqref{eq:splitting-direct-sum} that
	\[
	\sum_{j=1}^N|\lambda_j(H_{V_n})|^\gamma
	=
	\sum_{j=1}^{k}|\lambda_j(H_{V_n^{\mathrm{in}}})|^\gamma
	+
	\sum_{j=1}^{N-k}|\lambda_j(H_{V_n^{\mathrm{out}}})|^\gamma
	+o_n(1),
	\]
	which completes the proof.
\end{proof}

\subsection{Existence of optimizers}\label{subsec:compactness}
With the spectral splitting result established above, we now prove the existence of an optimizer by analyzing normalized maximizing sequences.
\begin{proof}[Proof of Theorem~\ref{thm:existence}]
	Let $\{V_n\}$ be a maximizing sequence for the variational problem
	\eqref{eq:normalized-constant}. Then we have
	\begin{equation}\label{eq:maximizing}
		\|V_n\|_p^p=1,
		\qquad
		\sum_{j=1}^N|\lambda_j(H_{V_n})|^\gamma
		\to L^{(N)}_{\gamma,d,s}>0 .
	\end{equation}
	Since 
	$\sum_{j=1}^N|\lambda_j(H_{V_n})|^\gamma
	\leq N|\lambda_1(H_{V_n})|^\gamma,$
	\eqref{eq:maximizing} and \eqref{eq:nonvanishing} allow us, for $n$ large enough, to choose suitable translations such that
	\begin{equation}\label{eq:nonvanishing-mass}
		\int_{B_1}V_n^p\,dx\geq c>0 .
	\end{equation}
	Define
	\[
	\alpha:=
	\lim_{R\to\infty}\liminf_{n\to\infty}
	\int_{B_R}V_n^p\,dx .
	\]
	Using \eqref{eq:maximizing} and \eqref{eq:nonvanishing-mass},
	we have $\alpha\in(0,1]$, so vanishing cannot occur.
	
	If $\alpha=1$, then, by \eqref{eq:maximizing},
	\[
	\lim_{R\to\infty}\limsup_{n\to\infty}
	\int_{|x|>R}V_n^p\,dx
	=
	\lim_{R\to\infty}\limsup_{n\to\infty}
	\left(
	1-\int_{B_R}V_n^p\,dx
	\right) 
	=
	1-
	\lim_{R\to\infty}\liminf_{n\to\infty}
	\int_{B_R}V_n^p\,dx 
	=0.
	\]	
	Therefore, Corollary~\ref{cor:tight-maximizing} implies that, up to a subsequence, $\{V_n\}$ converges strongly in $L^p(\mathbb R^d)$ to an optimizer, thereby proving the existence of optimizers for $L^{(N)}_{\gamma,d,s}$.
	
	Assume now that $0<\alpha<1$. 
	Passing to a subsequence and using a diagonal argument, we may choose $R_n\to\infty$ such that
	\begin{equation}\label{eq:measure-split}
		\|V_n^{\mathrm{in}}\|_p^p=\int_{B_{R_n}}V_n^p\,dx\to\alpha,
		\qquad
		\int_{B_{2R_n}\setminus B_{R_n}}V_n^p\,dx\to0 .
	\end{equation}
	Applying Lemma~\ref{lem:splitting} together with
	\eqref{eq:maximizing} and \eqref{eq:measure-split}, we obtain
	\begin{equation*}
		\begin{aligned}
			L^{(N)}_{\gamma,d,s}
			&=
			\lim_{n\to\infty}
			\left(
			\sum_{j=1}^k|\lambda_j(H_{V_n^{\mathrm{in}}})|^\gamma
			+
			\sum_{j=1}^{N-k}
			|\lambda_j(H_{V_n^{\mathrm{out}}})|^\gamma
			\right)\\
			&\leq
			\lim_{n\to\infty}
			\left(
			L^{(k)}_{\gamma,d,s}\|V_n^{\mathrm{in}}\|_p^p
			+
			L^{(N-k)}_{\gamma,d,s}
			\|V_n^{\mathrm{out}}\|_p^p
			\right)\\
			&=
			\alpha L^{(k)}_{\gamma,d,s}
			+(1-\alpha)L^{(N-k)}_{\gamma,d,s}
			\leq
			L^{(N)}_{\gamma,d,s}.
		\end{aligned}
	\end{equation*}
	Therefore, all inequalities above are equalities.  Since $0<\alpha<1$
	and $L^{(N)}_{\gamma,d,s}>0$, the cases $k=0$ and $k=N$ are excluded. Accordingly, we deduce that $1\leq k\leq N-1$, \begin{equation}\label{eq:inner-maximizing}
		L^{(k)}_{\gamma,d,s}\|V_n^{\mathrm{in}}\|_p^p
		-\sum_{j=1}^k|\lambda_j(H_{V_n^{\mathrm{in}}})|^\gamma
		\to0,
	\end{equation}
	and 
	\begin{equation}\label{eq:plateau-split}
		L^{(k)}_{\gamma,d,s}
		=
		L^{(N-k)}_{\gamma,d,s}
		=
		L^{(N)}_{\gamma,d,s}.
	\end{equation}
	We next verify that $\{V_n^{\mathrm{in}}\}$ is tight. For every fixed $R>0$, the convergence $R_n\to\infty$ ensures that $R_n>R$ for all sufficiently large $n$, so
	\[
	\int_{|x|>R}(V_n^{\mathrm{in}})^p\,dx
	=
	\int_{B_{R_n}\setminus B_R}V_n^p\,dx 
	=
	\int_{B_{R_n}}V_n^p\,dx
	-
	\int_{B_R}V_n^p\,dx.
	\]	
	Using \eqref{eq:measure-split}, we obtain
	\[
	\limsup_{n\to\infty}
	\int_{|x|>R}(V_n^{\mathrm{in}})^p\,dx
	=
	\alpha-
	\liminf_{n\to\infty}
	\int_{B_R}V_n^p\,dx.
	\]
	Letting $R\to\infty$ and recalling the definition of $\alpha$, we find
	\[
	\lim_{R\to\infty}\limsup_{n\to\infty}
	\int_{|x|>R}(V_n^{\mathrm{in}})^p\,dx
	=0.
	\]	
	Thus, \eqref{eq:tight} holds for $\{V_n^{\mathrm{in}}\}$.
	From \eqref{eq:measure-split}, \eqref{eq:inner-maximizing}, and
	Corollary~\ref{cor:tight-maximizing}, it then follows that, up to a subsequence, $V_n^{\mathrm{in}}\to V$ strongly in $L^p(\mathbb R^d)$, where $V$ is an optimizer for $L^{(k)}_{\gamma,d,s}$.	
	Observe that, by \eqref{eq:plateau-split},
	\[
	L^{(N)}_{\gamma,d,s}\|V\|_p^p
	=
	L^{(k)}_{\gamma,d,s}\|V\|_p^p
	=
	\sum_{j=1}^k|\lambda_j(H_V)|^\gamma
	\leq
	\sum_{j=1}^N|\lambda_j(H_V)|^\gamma
	\leq
	L^{(N)}_{\gamma,d,s}\|V\|_p^p,
	\]
	which shows that $V$ attains $L_{\gamma,d,s}^{(N)}$
	in the quotient formulation \eqref{eq:constant}.
	By scaling invariance, we obtain a normalized optimizer,
	completing the proof.
\end{proof}

\section{Euler--Lagrange equation and asymptotic behavior at infinity}\label{sec:structure}
This section is devoted to proving Theorem~\ref{thm:structure}.
In Subsection~\ref{subsec:EL}, we establish a resolvent representation and prove positivity almost everywhere and the spectral gap, which are essential for deriving the Euler--Lagrange equation.
In Subsection~\ref{subsec:tail}, we investigate the regularity
of the optimizer and the corresponding eigenfunctions,
and derive their precise decay estimates.

\subsection{Positivity, spectral gap, and Euler--Lagrange equation}\label{subsec:EL}
We begin with the following resolvent representation for
eigenfunctions associated with negative eigenvalues.

\begin{lem}\label{lem:resolvent-representation}
	Assume that $V\in L^p(\mathbb R^d)$ is nonnegative. Let
	$u\in H^s(\mathbb R^d)$ be real-valued and satisfy
	$H_Vu=-au$ in the weak sense for some $a>0$. Then $Vu\in L^{\frac{2p}{p+1}}(\mathbb R^d)
	\cap H^{-s}(\mathbb R^d)$, and
	\begin{equation}\label{eq:resolvent-equation}
		u=G_a*(Vu)
	\end{equation}
	almost everywhere in $\mathbb R^d$, where $G_a$ is the convolution
	kernel of $\left((-\Delta)^s+a\right)^{-1}$ defined in
	\eqref{eq:kernel-coefficient}.	
\end{lem}
\begin{proof}
	First, it follows from H\"older's inequality and \eqref{eq:GN} that
	\begin{equation*}\label{eq:Vu-Lq}
		\|Vu\|_{\frac{2p}{p+1}}
		\leq
		\|V\|_p
		\|u\|_{\frac{2p}{p-1}}
		\leq
		C\|V\|_p\|u\|_{H^s}.
	\end{equation*}
	Moreover, for every $\varphi\in H^s(\mathbb R^d)$, we have
	\[
	\left|\int_{\mathbb R^d}Vu\varphi\,dx\right|
	\leq
	\|V\|_p
	\|u\varphi\|_{\frac{p}{p-1}}
	\leq
	\|V\|_p
	\|u\|_{\frac{2p}{p-1}}
	\|\varphi\|_{\frac{2p}{p-1}}
	\leq
	C\|V\|_p\|u\|_{H^s}\|\varphi\|_{H^s}.	
	\]
	Hence $Vu\in L^{\frac{2p}{p+1}}(\mathbb R^d)\cap H^{-s}(\mathbb R^d)$.
	
	It remains to prove \eqref{eq:resolvent-equation}. Choose $V_m,u_m\in C_c^\infty(\mathbb R^d)$ such that
	\[
	V_m\to V\quad\text{in }L^p(\mathbb R^d),
	\qquad
	u_m\to u\quad\text{in }H^s(\mathbb R^d).
	\]
	By H\"older's inequality and the continuous embedding
	$H^s(\mathbb R^d)\hookrightarrow
	L^{\frac{2p}{p-1}}(\mathbb R^d)$, we obtain
	\begin{equation}\label{eq:product-Lq-convergence}
		\|V_mu_m-Vu\|_{\frac{2p}{p+1}}
		\leq
		\|V_m-V\|_p
		\|u_m\|_{\frac{2p}{p-1}}
		+
		\|V\|_p
		\|u_m-u\|_{\frac{2p}{p-1}}
		\to0.
	\end{equation}
	Furthermore, for every $\varphi\in H^s(\mathbb R^d)$ with
	$\|\varphi\|_{H^s}\leq1$, it holds that
	\[
	\left|
	\int_{\mathbb R^d}(V_mu_m-Vu)\varphi\,dx
	\right|
	\leq
	\|V_m-V\|_p
	\|u_m\|_{\frac{2p}{p-1}}
	\|\varphi\|_{\frac{2p}{p-1}}  
	+
	\|V\|_p
	\|u_m-u\|_{\frac{2p}{p-1}}
	\|\varphi\|_{\frac{2p}{p-1}}
	\to0,
	\]
	which, by the definition of the $H^{-s}$ norm, yields
	\begin{equation}\label{eq:product-Hminus-s}
		\|V_mu_m-Vu\|_{H^{-s}}\to0.
	\end{equation}
	
	Set $f_m:=V_mu_m$. In view of the definition of $G_a$, one has
	\begin{equation}\label{eq:resolvent-convolution-fm}
		\left((-\Delta)^s+a\right)^{-1}f_m=G_a*f_m.
	\end{equation}
	Define 
	\[w_m:=\left((-\Delta)^s+a\right)^{-1}(f_m-Vu).\]
	Then $\left((-\Delta)^s+a\right)w_m=f_m-Vu$ in the weak sense.
	Testing this equation with $w_m$ implies
	\[
	\|(-\Delta)^{s/2}w_m\|_2^2+a\|w_m\|_2^2
	=
	\langle f_m-Vu,w_m\rangle_{H^{-s},H^s}
	\leq
	\|f_m-Vu\|_{H^{-s}}\|w_m\|_{H^s}.
	\]
	Observe that
	$\|(-\Delta)^{s/2}w_m\|_2^2+a\|w_m\|_2^2
	\geq\min\{1,a\}\|w_m\|_{H^s}^2$.
	It thus follows from \eqref{eq:product-Hminus-s} that
	\[
	\|w_m\|_{H^s}
	\leq
	\frac{1}{\min\{1,a\}}
	\|f_m-Vu\|_{H^{-s}}
	\to0.
	\]
	Hence, by the definition of $w_m$,
	\begin{equation}\label{eq:resolvent-Hs-convergence}
		\left((-\Delta)^s+a\right)^{-1}f_m
		\to
		\left((-\Delta)^s+a\right)^{-1}(Vu)
		\quad\text{in }H^s(\mathbb R^d).
	\end{equation}
	On the other hand, since $G_a\in L^1(\mathbb R^d)$ (see \cite[Lemma~C.1(iv)]{FLS16}), we deduce from Young's inequality and \eqref{eq:product-Lq-convergence} that
	\begin{equation}\label{eq:convolution-Lq-convergence}
		\|G_a*f_m-G_a*(Vu)\|_{\frac{2p}{p+1}}
		\leq
		\|G_a\|_1
		\|f_m-Vu\|_{\frac{2p}{p+1}}
		\to 0.
	\end{equation}
	Using \eqref{eq:resolvent-Hs-convergence} and
	\eqref{eq:convolution-Lq-convergence}, both sides of
	\eqref{eq:resolvent-convolution-fm} converge in
	$\mathcal D'(\mathbb R^d)$. Therefore, by uniqueness of the distributional limit,
	\begin{equation}\label{eq:resolvent-convolution-limit}
		\left((-\Delta)^s+a\right)^{-1}(Vu)=G_a*(Vu)
	\end{equation}
	in $\mathcal D'(\mathbb R^d)$.
	Moreover, the eigenvalue equation
	$H_Vu=-au$ holds in the weak sense; that is,
	\[
	\big\langle (-\Delta)^{s/2}u,(-\Delta)^{s/2}\varphi\big\rangle
	+a\int_{\mathbb R^d}u\varphi\,dx
	=
	\int_{\mathbb R^d}Vu\varphi\,dx
	\qquad
	\text{for every }\varphi\in H^s(\mathbb R^d).
	\]
	Thus, the invertibility of $(-\Delta)^s+a$ from $H^s$ to $H^{-s}$ implies that
	\begin{equation*}
		u=\left((-\Delta)^s+a\right)^{-1}(Vu).
	\end{equation*}
	Combining this with \eqref{eq:resolvent-convolution-limit}, we obtain
	$u=G_a*(Vu)$
	almost everywhere, which completes the proof. 
\end{proof}

We now prove Theorem~\ref{thm:structure}-\ref{item:structure-EL},
concerning the spectral gap and the Euler--Lagrange equation.
Assume that $V\geq0$ is an optimizer for $L^{(N)}_{\gamma,d,s}$ with
$\|V\|_p^p=1$, and set
$K:=\#\{1\leq j\leq N:\lambda_j(H_V)<0\}$, where eigenvalues are
counted with multiplicity.

\begin{proof}[Proof of Theorem~\ref{thm:structure}-\ref{item:structure-EL}]
The proof is organized into the following three steps.

\textbf{Step 1. Positivity almost everywhere.}
Let $u$ be a real eigenfunction associated with the lowest eigenvalue $\lambda_1(H_V)$. 
Noting from \eqref{eq:kernel-form} that
$\|(-\Delta)^{s/2}|u|\|_2^2\leq\|(-\Delta)^{s/2}u\|_2^2$,
while the $L^2$-norm and the potential term are unchanged, the
min--max principle implies that $|u|$ is also an eigenfunction associated with $\lambda_1(H_V)$. Hence we may assume that $u\geq0$.

Write $\lambda_1(H_V)=-a$ with $a>0$.
By Lemma~\ref{lem:resolvent-representation}, we get
\[
u=G_a*(Vu).
\]
Since $G_a>0$ (see \cite[Lemma~C.1(i)]{FLS16}) and $Vu\geq0$ is not identically zero, the above identity implies that 
\[
u(x)=\int_{\mathbb R^d}G_a(x-y)V(y)u(y)\,dy>0
\]
for almost every $x\in\mathbb R^d$.

We next show that the lowest eigenvalue $\lambda_1(H_V)$ is simple. Suppose otherwise that there exists a real eigenfunction $v\not\equiv0$ associated with $\lambda_1(H_V)$ and orthogonal to $u$. Since $u>0$ almost everywhere and $\int_{\mathbb R^d}uv\,dx=0$, it follows that
\[
|\{x\in\mathbb R^d:v(x)>0\}|>0
\qquad\text{and}\qquad
|\{x\in\mathbb R^d:v(x)<0\}|>0,
\]
which, together with \eqref{eq:kernel-form}, means that
\[
\|(-\Delta)^{s/2}v\|_2^2
-\|(-\Delta)^{s/2}|v|\|_2^2
=
c_{d,s}
\iint_{\mathbb R^d\times\mathbb R^d}
\frac{|v(x)v(y)|-v(x)v(y)}
{|x-y|^{d+2s}}\,dx\,dy
>0.
\]
Consequently,
\[
q_V[|v|]
<
q_V[v]
=
\lambda_1(H_V)\|v\|_2^2.
\]
This contradicts the min--max characterization of
$\lambda_1(H_V)$.
Therefore $\lambda_1(H_V)$ is simple, and we denote its normalized
eigenfunction by $u_1$, chosen so that $u_1>0$ almost everywhere.

Finally, we prove that $V>0$ almost everywhere. If not, then
$|\{x\in\mathbb R^d:V(x)=0\}|>0$.
In that case, there exists $R>0$ such that $|\{x\in B_R:V(x)=0\}|>0$.
Set
\[
h_0:=\mathbf 1_{\{x\in B_R:V(x)=0\}}.
\]
Then $h_0$ is nonzero and has bounded support. 
For $t>0$, consider the perturbation $V+th_0$. 
Since $\lambda_1(H_V)$ is simple, the Hellmann--Feynman formula gives
\begin{equation}\label{eq:first-eigenvalue-variation}
	\left.\frac{d}{dt}\lambda_1(H_{V+th_0})\right|_{t=0}
	=
	\left\langle
	u_1,
	\left.\frac{d}{dt}H_{V+th_0}\right|_{t=0}u_1
	\right\rangle
	=
	-\int_{\mathbb R^d}h_0 u_1^2\,dx<0,
\end{equation}
where we have used the fact that $u_1>0$ almost everywhere.
Noting that $H_{V+th_0}=H_V-th_0\leq H_V$ in the form sense, the
min--max principle yields $\lambda_j(H_{V+th_0})\leq\lambda_j(H_V)$ for $1\leq j\leq N$.
Hence
\begin{equation}\label{eq:Riesz-monotonicity}
	|\lambda_j(H_{V+th_0})|^\gamma
	\geq
	|\lambda_j(H_V)|^\gamma,
	\qquad 1\leq j\leq N.
\end{equation}
Moreover, by Taylor expansion,
\begin{equation}\label{eq:first-Riesz-variation}
	|\lambda_1(H_{V+th_0})|^\gamma
	=
	|\lambda_1(H_V)|^\gamma
	+
	\gamma|\lambda_1(H_V)|^{\gamma-1}
	\left(\int_{\mathbb R^d}h_0 u_1^2\,dx\right)t
	+o(t).
\end{equation}
On the other hand, using $Vh_0=0$ almost everywhere and $\|V\|_p^p=1$, we have
\begin{equation}\label{eq:denominator-variation}
	\|V+th_0\|_p^p=1+t^p\|h_0\|_p^p=1+o(t),
\end{equation}
where the last equality follows from $p>1$. Thus, we infer from  \eqref{eq:Riesz-monotonicity}, \eqref{eq:first-Riesz-variation},  \eqref{eq:denominator-variation} and \eqref{eq:first-eigenvalue-variation} that
\[
\begin{aligned}
	\frac{\displaystyle\sum_{j=1}^N
		|\lambda_j(H_{V+th_0})|^\gamma}
	{\|V+th_0\|_p^p}
	&\geq
	\frac{
		\displaystyle\sum_{j=1}^N|\lambda_j(H_V)|^\gamma
		+
		\gamma|\lambda_1(H_V)|^{\gamma-1}
		\left(\int_{\mathbb R^d}h_0u_1^2\,dx\right)t
		+o(t)}
	{1+o(t)}
	\\
	&=
	L^{(N)}_{\gamma,d,s}
	+
	\gamma|\lambda_1(H_V)|^{\gamma-1}
	\left(\int_{\mathbb R^d}h_0u_1^2\,dx\right)t
	+o(t)
	\\
	&>
	L^{(N)}_{\gamma,d,s}
\end{aligned}
\]
for sufficiently small $t>0$. This contradicts the definition of $L^{(N)}_{\gamma,d,s}$. So $V>0$ almost everywhere.

\textbf{Step 2. Spectral gap.}
If $K<N$, then $\lambda_K(H_V)<0=\lambda_{K+1}(H_V)$,
and hence \eqref{eq:gap} follows immediately.

It remains to consider $K=N$. 
Suppose, for a contradiction, that
$\lambda_K(H_V)=\lambda_{K+1}(H_V)=\lambda<0$.
Let $m_0\geq2$ be the multiplicity of $\lambda$, and let
$k_0$ denote the number of occurrences of $\lambda$ among the first
$K$ eigenvalues. Apparently, we have $1\leq k_0<m_0$.
Using the positivity of the optimizer established in Step~1, we define, for $m\in\mathbb N$,
\begin{equation}\label{eq:variation-sets}
	E_m:=\{x\in\mathbb R^d:|x|\leq m,\ V(x)\geq1/m\}.
\end{equation}
Then, for any real-valued bounded measurable function $h$ supported in $E_m$, it holds that
$V+th\geq0$
for all sufficiently small $|t|$.
Choose a real orthonormal basis
$v_1,\ldots,v_{m_0}$ of $\ker(H_V-\lambda)$, and let
\[
X_h:=
\left(\int_{\mathbb R^d}h v_a v_b\,dx\right)_{a,b=1}^{m_0}.
\]
The matrix $X_h$ is symmetric, and we denote its eigenvalues by
$\mu_1\leq\cdots\leq\mu_{m_0}$.
By perturbation theory (see \cite[Chapter~VIII, Section~2, Theorem~2.6]{Kat76}), the perturbed eigenvalues of $H_{V+th}$ near $\lambda$ admit the asymptotic expansions
\[
\lambda_a(t)=\lambda-t\mu_a+o(|t|),
\qquad 1\leq a\leq m_0.	
\]
Since $\lambda<0$, Taylor expansion yields
\begin{equation}\label{eq:eigenvalue-perturbation}
|\lambda_a(t)|^\gamma
=
|\lambda|^\gamma
+\gamma t|\lambda|^{\gamma-1}\mu_a
+o(|t|),
\qquad 1\leq a\leq m_0.
\end{equation}
For $t>0$, the $k_0$ smallest eigenvalues among
$\lambda_1(t),\ldots,\lambda_{m_0}(t)$ correspond, at first order,
to the largest $k_0$ eigenvalues of $X_h$.
For $t<0$, they correspond to the smallest $k_0$ eigenvalues of $X_h$.

Fix real orthonormal eigenfunctions $u_1,\ldots,u_{K-k_0}$ of $H_V$ corresponding to $\lambda_1(H_V),\ldots,\lambda_{K-k_0}(H_V)$, respectively. 
Since all eigenvalues strictly below $\lambda$ are included among
the first $K$ levels, we have
\begin{equation}\label{eq:lower-spectral-expansion}
	\sum_{j=1}^{K-k_0}|\lambda_j(H_{V+th})|^\gamma
	=
	\sum_{j=1}^{K-k_0}|\lambda_j(H_V)|^\gamma
	+
	\gamma t\sum_{j=1}^{K-k_0}
	|\lambda_j(H_V)|^{\gamma-1}
	\int_{\mathbb R^d}hu_j^2\,dx
	+o(|t|).
\end{equation}
The coefficient of $t$ is independent of the choice of orthonormal
eigenfunctions within each eigenspace, by the invariance of the trace.
Moreover, Taylor expansion gives
\begin{equation}\label{eq:denominator-first-order}
	\|V+th\|_p^p
	=
	1+pt\int_{\mathbb R^d}V^{p-1}h\,dx+o(|t|),
\end{equation}
and
\begin{equation}\label{eq:denominator-reciprocal-expansion}
	\left(1+pt\int_{\mathbb R^d}V^{p-1}h\,dx+o(t)\right)^{-1}
	=
	1-pt\int_{\mathbb R^d}V^{p-1}h\,dx+o(t).
\end{equation}
Set
\[
c_h:=
\frac{
	pL^{(N)}_{\gamma,d,s}
	\displaystyle\int_{\mathbb R^d}V^{p-1}h\,dx
	-
	\gamma\displaystyle\sum_{j=1}^{K-k_0}
	|\lambda_j(H_V)|^{\gamma-1}
	\int_{\mathbb R^d}h u_j^2\,dx
}{
	\gamma|\lambda|^{\gamma-1}
}.
\]
For $t>0$ small enough, using \eqref{eq:lower-spectral-expansion},
\eqref{eq:denominator-first-order},  \eqref{eq:eigenvalue-perturbation}, and \eqref{eq:denominator-reciprocal-expansion}, we obtain
\[
\begin{aligned}
	&\frac{\displaystyle\sum\limits_{j=1}^K
		|\lambda_j(H_{V+th})|^\gamma}
	{\|V+th\|_p^p}\\
	&=
	\frac{
		L^{(N)}_{\gamma,d,s}
		+t\left(
		\gamma\displaystyle\sum\limits_{j=1}^{K-k_0}
		|\lambda_j(H_V)|^{\gamma-1}
		\int_{\mathbb R^d}hu_j^2\,dx
		+
		\gamma|\lambda|^{\gamma-1}
		\displaystyle\sum\limits_{a=m_0-k_0+1}^{m_0}\mu_a
		\right)
		+o(t)}
	{
		1+pt\displaystyle\int_{\mathbb R^d}V^{p-1}h\,dx+o(t)
	}\\
	&=
	L^{(N)}_{\gamma,d,s}
	+t\left(
	\gamma\sum_{j=1}^{K-k_0}
	|\lambda_j(H_V)|^{\gamma-1}
	\int_{\mathbb R^d}hu_j^2\,dx
	+
	\gamma|\lambda|^{\gamma-1}
	\sum_{a=m_0-k_0+1}^{m_0}\mu_a
	-
	pL^{(N)}_{\gamma,d,s}
	\int_{\mathbb R^d}V^{p-1}h\,dx
	\right)
	+o(t)\\
	&=
	L^{(N)}_{\gamma,d,s}
	+
	\gamma|\lambda|^{\gamma-1}
	\left(
	\sum_{a=m_0-k_0+1}^{m_0}\mu_a-c_h
	\right)t
	+o(t),
\end{aligned}
\]
which, together with the fact that $V$ attains $L^{(N)}_{\gamma,d,s}$, implies that
\begin{equation}\label{eq:right-slope}
	\gamma|\lambda|^{\gamma-1}
	\left(
	\sum_{a=m_0-k_0+1}^{m_0}\mu_a-c_h
	\right)\leq0.
\end{equation}
Furthermore, for sufficiently small $t<0$, the same expansions yield
\[
\begin{aligned}
	\frac{\displaystyle\sum\limits_{j=1}^K
		|\lambda_j(H_{V+th})|^\gamma}
	{\|V+th\|_p^p}
	&=
	\frac{
		L^{(N)}_{\gamma,d,s}
		+t\left(
		\gamma\displaystyle\sum\limits_{j=1}^{K-k_0}
		|\lambda_j(H_V)|^{\gamma-1}
		\int_{\mathbb R^d}hu_j^2\,dx
		+
		\gamma|\lambda|^{\gamma-1}
		\displaystyle\sum\limits_{a=1}^{k_0}\mu_a
		\right)
		+o(|t|)}
	{
		1+pt\displaystyle\int_{\mathbb R^d}V^{p-1}h\,dx
		+o(|t|)
	}\\
	&=
	L^{(N)}_{\gamma,d,s}
	+
	\gamma|\lambda|^{\gamma-1}
	\left(
	\sum_{a=1}^{k_0}\mu_a-c_h
	\right)t
	+o(|t|),
\end{aligned}
\]
which means that
\begin{equation}\label{eq:left-slope}
	\gamma|\lambda|^{\gamma-1}
	\left(
	\sum_{a=1}^{k_0}\mu_a-c_h
	\right)\geq0.
\end{equation}
Noting that $\gamma|\lambda|^{\gamma-1}>0$, we conclude from \eqref{eq:right-slope} and \eqref{eq:left-slope} that
\begin{equation}\label{eq:slope}
	\sum_{a=m_0-k_0+1}^{m_0}\mu_a
	\leq c_h
	\leq
	\sum_{a=1}^{k_0}\mu_a.
\end{equation}
On the other hand, the ordering of the eigenvalues gives
$\sum_{a=1}^{k_0}\mu_a
\leq
\sum_{a=m_0-k_0+1}^{m_0}\mu_a.$
Combining this with \eqref{eq:slope}, we get
$0=\sum_{a=1}^{k_0}\left(\mu_{m_0-k_0+a}-\mu_a\right).$
As $0<k_0<m_0$, this forces $\mu_1=\cdots=\mu_{m_0}$. Hence $X_h$ is a scalar matrix for every such $h$.
In particular, for $a\neq b$,
\[
\int_{\mathbb R^d}h v_a v_b\,dx=0,
\qquad
\int_{\mathbb R^d}h(v_a^2-v_b^2)\,dx=0.
\]
Since $h$ is arbitrary among bounded measurable functions supported
in $E_m$, it follows that
\[
v_av_b=0,
\qquad
v_a^2=v_b^2
\quad\text{almost everywhere on }E_m.
\]
Moreover, by the definition of $E_m$ in
\eqref{eq:variation-sets}, we have
$E_m\subset E_{m+1}$ for every $m\in\mathbb N$.
In view of $V>0$ almost everywhere, this yields
	\begin{equation}\label{eq:variation-exhaustion}
		\left|\mathbb R^d\setminus\bigcup_{m=1}^\infty E_m\right|=0.
\end{equation}
So the above identities hold almost
everywhere in $\mathbb R^d$. Consequently,
$v_a^4=v_a^2v_b^2=(v_av_b)^2=0$ almost everywhere in $\mathbb R^d.$
Namely, $v_a=v_b=0$ almost everywhere, contradicting their
normalization. Therefore,
$\lambda_K(H_V)<\lambda_{K+1}(H_V)$,
which proves \eqref{eq:gap}.

\textbf{Step 3. Euler--Lagrange equation.}
By \eqref{eq:gap}, we know that the first $K$ eigenvalues are
strictly negative, and each eigenvalue among them is counted
with its full multiplicity.
Since $V$ attains $L^{(N)}_{\gamma,d,s}$ with $\|V\|_p^p=1$,
we have
\begin{equation}\label{eq:effective-rank}
	L^{(N)}_{\gamma,d,s}
=
\sum_{j=1}^N|\lambda_j|^\gamma
=
\sum_{j=1}^K|\lambda_j|^\gamma
\leq L^{(K)}_{\gamma,d,s}
\leq L^{(N)}_{\gamma,d,s}.
\end{equation}
Thus $V$ also attains $L^{(K)}_{\gamma,d,s}$.
Denote by $u_1,\ldots,u_K$ a family of real orthonormal
eigenfunctions of $H_V$ corresponding to
$\lambda_1(H_V),\ldots,\lambda_K(H_V)$, respectively.
As in Step~2, let $h$ be any real-valued bounded measurable
function supported in $E_m$.
By perturbation theory and \eqref{eq:effective-rank}, we obtain
\[
\begin{aligned}
	0
	=
	\left.
	\frac{d}{dt}
	\frac{\displaystyle\sum_{j=1}^K
		|\lambda_j(H_{V+th})|^\gamma}
	{\|V+th\|_p^p}
	\right|_{t=0}
	&=
	\gamma\sum_{j=1}^K
	|\lambda_j(H_V)|^{\gamma-1}
	\int_{\mathbb R^d}hu_j^2\,dx
	-
	p\left(\sum_{j=1}^K|\lambda_j(H_V)|^\gamma\right)
	\int_{\mathbb R^d}V^{p-1}h\,dx\\
	&=
	\gamma\sum_{j=1}^K
	|\lambda_j(H_V)|^{\gamma-1}
	\int_{\mathbb R^d}hu_j^2\,dx
	-
	pL^{(N)}_{\gamma,d,s}
	\int_{\mathbb R^d}V^{p-1}h\,dx\\
	&=
	\int_{\mathbb R^d}
	\left(
	\gamma\sum_{j=1}^K
	|\lambda_j(H_V)|^{\gamma-1}u_j^2
	-
	pL^{(N)}_{\gamma,d,s}V^{p-1}
	\right)h\,dx.
\end{aligned}
\]
Since $h$ is arbitrary among bounded measurable functions supported
in $E_m$, it follows that
\[
\gamma\sum_{j=1}^K
|\lambda_j(H_V)|^{\gamma-1}u_j^2
=
pL^{(N)}_{\gamma,d,s}V^{p-1}
\quad\text{almost everywhere on }E_m.
\]
Using \eqref{eq:variation-exhaustion}, we conclude that this identity holds almost everywhere in $\mathbb R^d$.
Solving for $V$ yields \eqref{eq:EL}.
Substituting this expression into $H_Vu_j=\lambda_j u_j$
gives \eqref{eq:EL-equation} in the weak sense. 
This completes the proof of Theorem~\ref{thm:structure}-\ref{item:structure-EL}.
\end{proof}

\subsection{Regularity and asymptotic behavior at infinity}\label{subsec:tail}
In this subsection, we establish the regularity and decay estimates of the optimizer $V$ and its associated eigenfunctions.
Recall that $u_1,\ldots,u_K$ are real orthonormal eigenfunctions
of $H_V$ corresponding to $\lambda_1,\ldots,\lambda_K$,
respectively, where $\lambda_j=\lambda_j(H_V)<0$ for
$1\leq j\leq K$.

\begin{proof}[Proof of Theorem~\ref{thm:structure}-\ref{item:structure-tail}]
We divide the proof into the following two steps.
	
\textbf{Step 1. Regularity and positivity.}
For $1\leq j\leq K$, set $a_j:=|\lambda_j|>0$.
Since $H_Vu_j=-a_j u_j$,
Lemma~\ref{lem:resolvent-representation} gives
\begin{equation}\label{eq:occupied-resolvent}
	u_j=G_{a_j}*(Vu_j)
	\quad\text{almost everywhere in }\mathbb R^d.
\end{equation}
Moreover, we recall from \cite[Lemma~C.1(iv)]{FLS16}
that, for every $a>0$,
\begin{equation}\label{eq:resolvent-integrability}
	G_a\in L^r(\mathbb R^d),
	\qquad
	1\leq r\leq\infty,\quad
	1-\frac1r<\frac{2s}{d}.
\end{equation}

We first establish the boundedness of each $u_j$ by an iteration
argument and then deduce the boundedness of $V$ from \eqref{eq:EL}.
By \eqref{eq:GN}, we have
$u_j\in L^{2p'}(\mathbb R^d)$ with $p'=p/(p-1)$. If $d\leq2s$, set $r_0=2p'$ and $1/t_0=1/p+1/r_0$.
Recalling that $p>\max\{1,d/(2s)\}$, we have
$\frac1{t_0}=\frac{p+1}{2p}<1\leq\frac{2s}{d}$.
Together with \eqref{eq:resolvent-integrability}, this yields
$G_{a_j}\in L^{t_0'}(\mathbb R^d)$
with $t_0'=t_0/(t_0-1)$.
It thus follows from \eqref{eq:occupied-resolvent}, H\"older's inequality, and Young's inequality that
\[
\|u_j\|_\infty
\leq\|G_{a_j}\|_{t_0'}\|Vu_j\|_{t_0}
\leq\|G_{a_j}\|_{t_0'}\|V\|_p\|u_j\|_{r_0}
<\infty.	
\]
On the other hand, suppose that $d>2s$.
In view of $p>\max\{1,d/(2s)\}$, set $\delta:=\frac{2s}{d}-\frac1p>0.$
We claim that, starting from $r_0=2p'$, there exists
some $n\in\mathbb N_0$ such that
\begin{equation}\label{eq:iteration-integrability}
	u_j\in L^{r_n}(\mathbb R^d),
\qquad \frac1{r_n}<\delta.
\end{equation}
Indeed, if $1/r_0<\delta$, the claim holds with $n=0$, since $u_j\in L^{r_0}(\mathbb R^d)$. Otherwise, for each $n\in\mathbb N_0$ such that $1/r_n\geq\delta$,
let $1/r_{n+1}=1/r_n-\delta/2$.
Choose $b>1$ such that
$1-\frac1b
=\frac1p+\frac{\delta}{2}
=\frac{2s}{d}-\frac{\delta}{2}
<\frac{2s}{d}$.
Then $G_{a_j}\in L^b(\mathbb R^d)$ by
\eqref{eq:resolvent-integrability}.
For each $n$ for which $r_n$ is defined, take $t_n>1$ such that
$1/t_n=1/p+1/r_n$.
Observe that
\[
	\frac1b+\frac1{t_n}
	=\left(1-\frac1p-\frac{\delta}{2}\right)
	+\left(\frac1p+\frac1{r_n}\right)
	=1-\frac{\delta}{2}+\frac1{r_n}
	=1+\frac1{r_{n+1}}.
\]
We deduce from \eqref{eq:occupied-resolvent},
Young's inequality, and H\"older's inequality that
\[
\|u_j\|_{r_{n+1}}
\leq\|G_{a_j}\|_b\|Vu_j\|_{t_n}
\leq\|G_{a_j}\|_b\|V\|_p\|u_j\|_{r_n}
<\infty,
\]
which, together with the fact that $u_j\in L^{r_0}(\mathbb R^d)$,
implies by induction that
$u_j\in L^{r_n}(\mathbb R^d)$.
Moreover, whenever $1/r_n\geq\delta$, we have
$1/r_{n+1}=1/r_n-\delta/2\geq\delta/2>0$.
Since $1/r_n=1/r_0-n\delta/2$ as long as the iteration
continues, it terminates after finitely many steps.
This proves the claim \eqref{eq:iteration-integrability}.
In view of \eqref{eq:iteration-integrability}, we have $0<\frac1{t_n}
=\frac1p+\frac1{r_n}
<\frac1p+\delta
=\frac{2s}{d}<1$.
Repeating the argument for the case $d\leq2s$, we obtain
\[
	\|u_j\|_\infty
	\leq\|G_{a_j}\|_{t_n'}\|Vu_j\|_{t_n}
	\leq\|G_{a_j}\|_{t_n'}\|V\|_p\|u_j\|_{r_n}
	<\infty,
\]
where $t_n'=t_n/(t_n-1)$. Consequently,
\[
u_j\in L^\infty(\mathbb R^d),
\qquad 1\leq j\leq K.
\]
Since $K$ is finite, \eqref{eq:EL} further implies that
$V\in L^\infty(\mathbb R^d)$.

We next establish further regularity of the eigenfunctions
and show that both $u_j$ and $V$ vanish at infinity.
Since $V\in L^\infty(\mathbb R^d)$ and
$u_j\in L^2(\mathbb R^d)$, the eigenvalue equation gives
\[
(-\Delta)^su_j=(V-a_j)u_j\in L^2(\mathbb R^d).
\]
By Plancherel's theorem,  we conclude that $u_j\in H^{2s}(\mathbb R^d)$.
Moreover, applying \cite[Proposition~B.3]{FLS16} to the bounded
potential $-V$, we have
\begin{equation}\label{eq:holder-regularity}
	u_j\in C^{0,\beta}(\mathbb R^d)
\qquad\text{for every }
0<\beta<\min\{1,2s\}.
\end{equation}
In particular, each $u_j$ is uniformly continuous.
We then claim that
$$u_j(x)\to0 \qquad\text{as } |x|\to\infty.$$
Indeed, otherwise there would exist $\varepsilon>0$ and
a sequence $|x_n|\to\infty$ such that
$|u_j(x_n)|\geq\varepsilon$.
By uniform continuity, there exists $\rho>0$,
independent of $n$, such that
\[
|u_j(x)-u_j(x_n)|<\frac{\varepsilon}{2}
\qquad\text{for all }x\in B_\rho(x_n).
\]
Consequently,
\[
|u_j(x)|
\geq |u_j(x_n)|-|u_j(x)-u_j(x_n)|
\geq\frac{\varepsilon}{2}
\qquad\text{for all }x\in B_\rho(x_n).
\]
Since $|x_n|\to\infty$, passing to a subsequence,
we may assume that the balls $B_\rho(x_n)$ are pairwise
disjoint. Consequently,
\[
\|u_j\|_2^2
\geq
\sum_{n=1}^{\infty}
\int_{B_\rho(x_n)}|u_j|^2\,dx
\geq
\sum_{n=1}^{\infty}
\frac{\varepsilon^2}{4}|B_\rho|
=\infty,
\]
a contradiction. Thus $u_j\in C_0(\mathbb R^d)$ for every $1\leq j\leq K$.
Since $u_j\in C_0(\mathbb R^d)$ for $1\leq j\leq K$
and $p>1$, \eqref{eq:EL} defines a continuous
representative of $V$ that vanishes at infinity.
Thus $V\in C_0(\mathbb R^d)$.

Finally, we show that $u_1(x)>0$ and $V(x)>0$
for every $x\in\mathbb R^d$.
Since $G_{a_j}\in L^1(\mathbb R^d)$ by \eqref{eq:resolvent-integrability} and
$Vu_j\in L^\infty(\mathbb R^d)$,
for every $z\in\mathbb R^d$, we have
\[
\begin{aligned}
	\sup_{x\in\mathbb R^d}
	\big|
	(G_{a_j}*(Vu_j))(x+z)
	-(G_{a_j}*(Vu_j))(x)
	\big|
	\leq
	\|Vu_j\|_\infty
	\|G_{a_j}(\cdot+z)-G_{a_j}\|_1
	\to0
	\text{ as }z\to0,
\end{aligned}
\]
where the convergence follows from the continuity
of translations in $L^1(\mathbb R^d)$.
Thus $G_{a_j}*(Vu_j)$ is continuous.
Since $u_j$ is also continuous by \eqref{eq:holder-regularity}
and \eqref{eq:occupied-resolvent} holds almost everywhere,
the identity holds everywhere. That is,
\[
u_j(x)=
\int_{\mathbb R^d}
G_{a_j}(x-y)V(y)u_j(y)\,dy,
\qquad x\in\mathbb R^d.
\]
As established in Step~1 of the proof of
Theorem~\ref{thm:structure}-\ref{item:structure-EL},
we have $u_1>0$ and $V>0$ almost everywhere.
In particular, $Vu_1\geq0$ is not identically zero.
Moreover, \cite[Lemma~C.1(i)]{FLS16} gives
$G_{a_1}(x-y)>0$ whenever $x\neq y$.
Consequently, for every $x\in\mathbb R^d$,
\[
u_1(x)=
\int_{\mathbb R^d}
G_{a_1}(x-y)V(y)u_1(y)\,dy>0,
\]
which, together with \eqref{eq:EL}, implies that
\[
	V(x)
	=\left(
	\frac{\gamma}{pL^{(N)}_{\gamma,d,s}}
	\sum_{j=1}^K|\lambda_j|^{\gamma-1}|u_j(x)|^2
	\right)^{1/(p-1)}
	\geq
	\left(
	\frac{\gamma}{pL^{(N)}_{\gamma,d,s}}
	|\lambda_1|^{\gamma-1}|u_1(x)|^2
	\right)^{1/(p-1)}
	>0,
	\qquad x\in\mathbb R^d.
\]

\textbf{Step 2. Asymptotic behavior at infinity.}
We first establish a decay estimate for each eigenfunction $u_j$.
Since $V\in C_0(\mathbb R^d)$, for every $1\leq j\leq K$,
there exists $R>0$ such that
\begin{equation}\label{eq:potential-exterior-bound}
	V(x)\leq\frac{a_j}{2},
\qquad |x|\geq R.
\end{equation}
We apply \cite[Lemma~C.2(i)]{FLS16} with potential $-V$,
eigenvalue $-a_j$, and parameter $\lambda=a_j/2$.
Indeed, $\|u_j\|_2=1$, $0<a_j/2<a_j$, and
\eqref{eq:potential-exterior-bound} implies that
$-V(x)+a_j/2\geq0$ for $|x|\geq R$.
It thus follows that, for some $C_j>0$,
\begin{equation}\label{eq:state-bound}
	|u_j(x)|
	\leq C_j(1+|x|)^{-(d+2s)},
	\qquad x\in\mathbb R^d.
\end{equation}

We next determine the precise asymptotic behavior at infinity
of each eigenfunction $u_j$.
Since $V\in L^\infty(\mathbb R^d)\cap C_0(\mathbb R^d)$,
we infer from \eqref{eq:state-bound} that
\[
	\|Vu_j\|_1
	\leq \|V\|_\infty\|u_j\|_1
	\leq C_j\|V\|_\infty
	\int_{\mathbb R^d}(1+|x|)^{-(d+2s)}\,dx
	\le C_j\|V\|_\infty C
	\int_0^\infty
	\frac{r^{d-1}}{(1+r)^{d+2s}}\,dr
	<\infty,
\]
and 
\[
|V(x)u_j(x)|
\leq C_j\|V\|_\infty
(1+|x|)^{-(d+2s)}.
\]
Moreover, we have
\[
	|x|^{d+2s}|V(x)u_j(x)|
	\leq C_j
	\left(\frac{|x|}{1+|x|}\right)^{d+2s}
	|V(x)|
	\leq C_j|V(x)|
	\to0
	\qquad\text{as }|x|\to\infty.
\]
For the resolvent kernel, \cite[Lemma~C.1(ii)--(iii)]{FLS16}
gives
\[
|G_{a_j}(x)|\leq C_j|x|^{-(d+2s)},
\qquad x\neq0.
\]
Indeed, part (ii) gives
$|G_{a_j}(x)|\leq C_j|x|^{-d}$ near the origin,
whereas part (iii) gives the required decay at infinity.
Moreover, by \eqref{eq:resolvent-integrability}
and \eqref{eq:kernel-coefficient}, we have
\[
G_{a_j}\in L^1(\mathbb R^d),
\qquad
\lim_{|x|\to\infty}|x|^{d+2s}G_{a_j}(x)
=a_{d,s}a_j^{-2},
\]
where the kernel limit follows from
\cite[Lemma~C.1(iii)]{FLS16}.
Thus the hypotheses of \cite[Lemma~C.3]{FLS16} are satisfied
with kernel $G_{a_j}$, function $Vu_j$, and exponent $\beta=d+2s>d$.
Applying that lemma to $u_j=G_{a_j}*(Vu_j)$, we obtain
\begin{equation}\label{eq:state-tail-limit}
	\lim_{|x|\to\infty}|x|^{d+2s}u_j(x)
=
a_{d,s}a_j^{-2}
\int_{\mathbb R^d}V(y)u_j(y)\,dy,
\end{equation}
which proves \eqref{eq:state-tail}.

Finally, we establish the precise asymptotic behavior of $V$ at infinity.
Since $K$ is finite, \eqref{eq:state-tail-limit} yields
\[
	\lim_{|x|\to\infty}
	|x|^{2(d+2s)}
	\sum_{j=1}^K|\lambda_j|^{\gamma-1}|u_j(x)|^2
	=
	\sum_{j=1}^K|\lambda_j|^{\gamma-1}
	\left(
	a_{d,s}a_j^{-2}
	\int_{\mathbb R^d}V(y)u_j(y)\,dy
	\right)^2.
\]
By \eqref{eq:EL}, we obtain, as $|x|\to\infty$,
\[
\begin{aligned}
	V(x)
	&=
	\left[
	\frac{\gamma}{pL^{(N)}_{\gamma,d,s}}
	\sum_{j=1}^K|\lambda_j|^{\gamma-1}
	\left(
	a_{d,s}a_j^{-2}
	\int_{\mathbb R^d}V(y)u_j(y)\,dy
	\right)^2
	\right]^{1/(p-1)}
	|x|^{-\frac{2(d+2s)}{p-1}}+
	o\left(|x|^{-\frac{2(d+2s)}{p-1}}\right).
\end{aligned}
\]
This completes the proof of Theorem~\ref{thm:structure}-\ref{item:structure-tail}.
\end{proof}

\section{Strict binding inequality}\label{sec:doubling}
This section is devoted to proving Theorem~\ref{thm:doubling}
by constructing a trial potential and using
a trace expansion.

\begin{proof}[Proof of Theorem~\ref{thm:doubling}]
Fix $N\in\mathbb N$ and denote by $V$ a normalized optimizer
for $L^{(N)}_{\gamma,d,s}$, whose existence follows from
Theorem~\ref{thm:existence}.
Let $K$ be as in Theorem~\ref{thm:structure}.
The argument proceeds in three steps.

\textbf{Step 1. The trial potential and interaction estimates.}
Recall the eigenpairs $(\lambda_j,u_j)_{j=1}^K$
from Theorem~\ref{thm:structure}, and set
\begin{equation}\label{eq:f-definition}
	f:=V^{p-1}
=\frac{\gamma}{pL^{(N)}_{\gamma,d,s}}
\sum_{j=1}^K|\lambda_j|^{\gamma-1}|u_j|^2.
\end{equation}
Clearly, $f\in C_0(\mathbb R^d)$, $f>0$ everywhere,
and $f\in L^{p/(p-1)}(\mathbb R^d)$ since $\|V\|_p^p=1$.
Moreover, integrating the expression for $f$ and using
$\|u_j\|_2=1$ for $1\leq j\leq K$, we get
\[
\|f\|_1
=
\frac{\gamma}{pL^{(N)}_{\gamma,d,s}}
\sum_{j=1}^K|\lambda_j|^{\gamma-1}
<\infty.
\]
Thus $f\in L^1(\mathbb R^d)\cap L^\infty(\mathbb R^d)$.

Let $e_1=(1,0,\ldots,0)\in\mathbb R^d$.
For $R>0$, define
\begin{equation}\label{eq:translated-functions}
	f_R(x):=f(x-Re_1),\qquad
u_{j,R}(x):=u_j(x-Re_1),\qquad 1\leq j\leq K,
\end{equation}
and the trial potential
\begin{equation}\label{eq:trial-potential}
	V_R(x):=(f(x)+f_R(x))^{\frac{1}{p-1}}
	=\left(V^{p-1}(x)+V^{p-1}(x-Re_1)\right)^{\frac{1}{p-1}}.
\end{equation}
Since $f,f_R\in L^{p/(p-1)}(\mathbb R^d)$, this yields $V_R\in L^p(\mathbb R^d)$.
We then introduce the interaction term
\begin{equation}\label{eq:interaction-term}
	A_R:=\|V_R\|_p^p-2\|V\|_p^p
=\int_{\mathbb R^d}
\left((f+f_R)^{\frac{p}{p-1}}
-f^{\frac{p}{p-1}}-f_R^{\frac{p}{p-1}}\right)\,dx,
\end{equation}
where we used translation invariance of the integral.
By using $\|V\|_p^p=1$, we obtain
\begin{equation}\label{eq:denominator}
	\|V_R\|_p^p=2+A_R.
\end{equation}
Furthermore, it follows that
\begin{equation}\label{eq:trial-uniform-bound}
	\|V_R\|_\infty
=\|(f+f_R)^{\frac{1}{p-1}}\|_\infty
\leq(\|f\|_\infty+\|f_R\|_\infty)^{\frac{1}{p-1}}
=2^{\frac{1}{p-1}}\|f\|_\infty^{\frac{1}{p-1}}
=2^{\frac{1}{p-1}}\|V\|_\infty,
\quad R>0,
\end{equation}	
so the trial potentials are uniformly bounded in $L^\infty(\mathbb R^d)$.
We now claim that
\begin{equation}\label{eq:overlap-convergence}
	A_R>0,\qquad A_R\to0\quad\text{as }R\to\infty.
\end{equation}
To this end, we first establish the elementary inequality
\begin{equation}\label{eq:overlap-upper}
	0\leq (a+b)^{\frac{p}{p-1}}-a^{\frac{p}{p-1}}-b^{\frac{p}{p-1}}
	\leq C_p\left(a^{\frac{1}{p-1}}b+ab^{\frac{1}{p-1}}\right),
	\qquad a,b\geq0.
\end{equation}
Indeed, when $a\geq b\geq0$, the mean value theorem yields
\[
	(a+b)^{\frac{p}{p-1}}
-a^{\frac{p}{p-1}}
-b^{\frac{p}{p-1}}
\leq
(a+b)^{\frac{p}{p-1}}
-a^{\frac{p}{p-1}}
\leq
\frac{p}{p-1}(a+b)^{\frac{1}{p-1}}b
\leq
\frac{p\,2^{\frac{1}{p-1}}}{p-1}
a^{\frac{1}{p-1}}b.
\]
The case $b\geq a$ follows by interchanging $a$ and $b$.
In particular, the first inequality in \eqref{eq:overlap-upper}
is strict whenever $a,b>0$.
Since $f,f_R>0$, we conclude that $A_R>0$.
To prove that $A_R\to0$ as $R\to\infty$,
we apply \eqref{eq:overlap-upper} to obtain
\[
\begin{aligned}
	A_R
	&\leq C_p\int_{\mathbb R^d}
	\left(f^{\frac{1}{p-1}}(x)f_R(x)
	+f(x)f_R^{\frac{1}{p-1}}(x)\right)\,dx\\
	&=C_p\int_{\mathbb R^d}
	\left(V(x)f_R(x)+f(x)V(x-Re_1)\right)\,dx\\
	&=C_p\int_{\mathbb R^d}
	f(x)\left(V(x+Re_1)+V(x-Re_1)\right)\,dx.
\end{aligned}
\]
For each fixed $x\in\mathbb R^d$,
$V(x+Re_1)+V(x-Re_1)\to0$ as $R\to\infty$,
since $V\in C_0(\mathbb R^d)$ by Theorem~\ref{thm:structure}.
Moreover, one can see that
\[
0\leq
f(x)\left(V(x+Re_1)+V(x-Re_1)\right)
\leq 2\|V\|_\infty f(x),
\]
where the right hand side belongs to $L^1(\mathbb R^d)$.
The dominated convergence theorem therefore implies that
$A_R\to0$ as $R\to\infty$.

We next show that
\begin{equation}\label{eq:overlap-small-o}
	R^{-2(d+2s)}=o(A_R)
\qquad \text{as } R\to\infty.
\end{equation}
First, for $0<b\leq a$, it follows from convexity and
the inequality $b^{p/(p-1)}\leq a^{1/(p-1)}b$ that
\begin{equation}\label{eq:scalar-overlap}
	(a+b)^{\frac{p}{p-1}}
	-a^{\frac{p}{p-1}}
	-b^{\frac{p}{p-1}}
	\geq
	\frac{1}{p-1}a^{\frac{1}{p-1}}b.
\end{equation}
Moreover, by Remark~\ref{rem:L1-threshold}, we have
\begin{equation*}
	f(x)\asymp |x|^{-2(d+2s)}
	\qquad \text{as } |x|\to\infty,
\end{equation*}
which implies that there exist $R_0\geq1$ and
$c_0,C_0>0$ such that
\begin{equation}\label{eq:f-tail-bounds}
	c_0|x|^{-2(d+2s)}
	\leq f(x)
	\leq C_0|x|^{-2(d+2s)},
	\qquad |x|\geq R_0.
\end{equation}
Choose a fixed $\eta\in(0,1/2)$ sufficiently small that
\[
\frac{C_0}{c_0}
\left(\frac{\eta}{1-\eta}\right)^{2(d+2s)}
\leq1.
\]
For $R>R_0/\eta>2R_0$ and $R_0\leq |x|\leq\eta R$,
the triangle and reverse triangle inequalities yield
\begin{equation}\label{eq:annulus-distance}
	R/2<(1-\eta)R\leq R-|x|
	\leq |x-Re_1|
	\leq R+|x|
	\leq(1+\eta)R.
\end{equation}
In particular, both $|x|$ and $|x-Re_1|$ are at least $R_0$.
Applying \eqref{eq:f-tail-bounds} and \eqref{eq:annulus-distance},
we infer that
\[
\frac{f_R(x)}{f(x)}
=\frac{f(x-Re_1)}{f(x)}
\leq
\frac{C_0}{c_0}
\left(\frac{|x|}{|x-Re_1|}\right)^{2(d+2s)}
\leq
\frac{C_0}{c_0}
\left(\frac{\eta}{1-\eta}\right)^{2(d+2s)}
\leq1.
\]
On the annulus
$\{x\in\mathbb R^d:R_0\leq |x|\leq\eta R\}$,
\eqref{eq:f-tail-bounds} and \eqref{eq:annulus-distance}
also give
\begin{equation}\label{eq:translated-tail-lower}
	f_R(x)\geq c_0|x-Re_1|^{-2(d+2s)}
\geq c_0(1+\eta)^{-2(d+2s)}R^{-2(d+2s)}>0.
\end{equation}
We may therefore apply \eqref{eq:scalar-overlap}
with $a=f(x)$ and $b=f_R(x)$ on this annulus.
Using \eqref{eq:f-tail-bounds} again, we find
\[
V(x)=f^{\frac{1}{p-1}}(x)
\geq
c_0^{\frac{1}{p-1}}
|x|^{-\frac{2(d+2s)}{p-1}},
\qquad |x|\geq R_0.
\]
Combining this bound with \eqref{eq:translated-tail-lower}
and integrating in polar coordinates, we obtain
\begin{equation}\label{eq:overlap-integral}
	\begin{aligned}
		A_R
		&\geq
		\frac{1}{p-1}
		\int_{\{R_0\leq |x|\leq\eta R\}}
		V(x)f_R(x)\,dx\\
		&\geq
		cR^{-2(d+2s)}
		\int_{\{R_0\leq |x|\leq\eta R\}}
		|x|^{-\frac{2(d+2s)}{p-1}}\,dx\\
		&\ge 
		cR^{-2(d+2s)}
		\int_{R_0}^{\eta R}
		r^{d-1-\frac{2(d+2s)}{p-1}}\,dr,
	\end{aligned}
\end{equation}
with constants independent of $R$.
The last integral in \eqref{eq:overlap-integral} can be evaluated explicitly:
\[
\int_{R_0}^{\eta R}
r^{d-1-\frac{2(d+2s)}{p-1}}\,dr
=
\begin{cases}
	\displaystyle
	\frac{
		(\eta R)^{d-\frac{2(d+2s)}{p-1}}
		-R_0^{d-\frac{2(d+2s)}{p-1}}
	}{
		d-\frac{2(d+2s)}{p-1}
	},
	& p>3+\dfrac{4s}{d},\\[3ex]
	\displaystyle
	\log\left(\frac{\eta R}{R_0}\right),
	& p=3+\dfrac{4s}{d}.
\end{cases}
\]
Observe that the condition $p\geq 3+\frac{4s}{d}$ in
\eqref{eq:doubling-range} is equivalent to
$d-\frac{2(d+2s)}{p-1}\geq0$.
Consequently, for all sufficiently large $R$,
\begin{equation*}
	A_R\geq
	\begin{cases}
		cR^{-2(d+2s)+d-\frac{2(d+2s)}{p-1}},
		& p>3+\dfrac{4s}{d},\\[2ex]
		cR^{-2(d+2s)}\log R,
		& p=3+\dfrac{4s}{d}.
	\end{cases}
\end{equation*}
In either case, $R^{2(d+2s)}A_R\to\infty$
as $R\to\infty$, which proves \eqref{eq:overlap-small-o}.

\textbf{Step 2. Orthogonalization and trace expansion.}
Consider the family
\[
(u_1,\ldots,u_K,u_{1,R},\ldots,u_{K,R}),
\]
where $u_{j,R}$ is defined in \eqref{eq:translated-functions}
and is an eigenfunction of $H_{V(\cdot-Re_1)}$
corresponding to the eigenvalue $\lambda_j$,
for $1\leq j\leq K$.
By \eqref{eq:state-bound} and $d+2s>d$, we readily see
that each eigenfunction belongs to $L^1(\mathbb R^d)$.
For $R\geq1$, we split the integral into
$\{|x|\geq R/2\}$ and $\{|x|<R/2\}$.
On the latter set, the reverse triangle inequality gives $|x-Re_1|\geq R-|x|>R/2$.
Hence, using \eqref{eq:state-bound} again, we obtain
\begin{equation}\label{eq:state-overlap}
	\begin{aligned}
		\int_{\mathbb R^d}|u_j(x)u_{k,R}(x)|\,dx
		&=
		\int_{\{|x|\geq R/2\}}
		|u_j(x)u_{k,R}(x)|\,dx+
		\int_{\{|x|<R/2\}}
		|u_j(x)u_{k,R}(x)|\,dx\\
		&\leq
		CR^{-(d+2s)}
		\bigl(\|u_k\|_1+\|u_j\|_1\bigr)\\
		&\leq CR^{-(d+2s)},
		\qquad 1\leq j,k\leq K.
	\end{aligned}
\end{equation}
Let $G_R$ be the Gram matrix of the ordered family above.
Since each of the two families is orthonormal, we have
\begin{equation}\label{eq:Gram}
G_R=I_{2K}+E_R,\qquad
E_R=
\begin{pmatrix}
	0 & (\langle u_j,u_{k,R}\rangle)_{j,k=1}^K\\
	(\langle u_{j,R},u_k\rangle)_{j,k=1}^K & 0
\end{pmatrix}.
\end{equation}
From \eqref{eq:state-overlap} and the fact that $K$ is fixed,
it follows that
\begin{equation}\label{eq:Gram-error-bound}
	\|E_R\|_{\mathrm{op}}
\leq 2K\max_{1\leq j,k\leq K}
|\langle u_j,u_{k,R}\rangle|
\leq CR^{-(d+2s)},
\end{equation}
which implies that $\|E_R\|_{\mathrm{op}}<1/2$ for all
sufficiently large $R$.
Noting that $E_R$ is symmetric, for every $\xi\in\mathbb R^{2K}$, one has
\[
\xi^{\mathsf T}G_R\xi
=|\xi|^2+\xi^{\mathsf T}E_R\xi
\geq(1-\|E_R\|_{\mathrm{op}})|\xi|^2
\geq\frac12|\xi|^2.
\]
Thus $G_R\geq\frac12 I_{2K}$, and so $G_R$
is positive definite.
Moreover, for any $a_j,b_j\in\mathbb R$, $1\leq j\leq K$,
we also have
\[
\begin{aligned}
		\left\|\sum_{j=1}^K(a_ju_j+b_ju_{j,R})\right\|_2^2
	&=\left\langle\sum_{j=1}^K(a_ju_j+b_ju_{j,R}),\sum_{k=1}^K(a_ku_k+b_ku_{k,R})\right\rangle\\
	&=\xi^{\mathsf T}G_R\xi
	\geq\frac12|\xi|^2
	=\frac12\sum_{j=1}^K(a_j^2+b_j^2),
\end{aligned}
\]
where $\xi=(a_1,\ldots,a_K,b_1,\ldots,b_K)^{\mathsf T}$.
Consequently, the family $(u_1,\ldots,u_K,u_{1,R},\ldots,u_{K,R})$ is linearly independent.

Recalling the definition of $V_R$ in
\eqref{eq:trial-potential}, we define the matrices
$P_R^0$ and $P_R^1$ by
\[
(P_R^0)_{jk}
:=\int_{\mathbb R^d}
(V_R(x)-V(x))u_j(x)u_k(x)\,dx
\]
and
\[
(P_R^1)_{jk}
:=\int_{\mathbb R^d}
(V_R(x)-V(x-Re_1))
u_{j,R}(x)u_{k,R}(x)\,dx,
\]
for $1\leq j,k\leq K$.
We then set
\begin{equation*}
	P_R=
	\begin{pmatrix}
		P_R^0 & 0\\
		0 & P_R^1
	\end{pmatrix}.
\end{equation*}
In view of \eqref{eq:trial-potential}, it follows that
$V_R\geq V$ and $V_R\geq V(\cdot-Re_1)$.
Therefore, for every
$\widehat{\xi}=(\widehat{\xi}_1,\ldots,\widehat{\xi}_K)^{\mathsf T}
\in\mathbb R^K$,
\[
\widehat{\xi}^{\mathsf T}P_R^0\widehat{\xi}
=\int_{\mathbb R^d}
(V_R(x)-V(x))
\left|\sum_{j=1}^K\widehat{\xi}_j u_j(x)\right|^2\,dx
\geq0,
\]
and
\[
\widehat{\xi}^{\mathsf T}P_R^1\widehat{\xi}
=\int_{\mathbb R^d}
(V_R(x)-V(x-Re_1))
\left|\sum_{j=1}^K\widehat{\xi}_j u_{j,R}(x)\right|^2\,dx
\geq0.
\]
Hence $P_R^0$, $P_R^1$, and $P_R$ are positive semidefinite, and all their eigenvalues are nonnegative.
Moreover, using \eqref{eq:f-definition} and its translated version,
we obtain
\begin{equation}\label{eq:trace-identity}
	\begin{aligned}
		&\frac{\gamma}{pL^{(N)}_{\gamma,d,s}}
		\sum_{j=1}^K|\lambda_j|^{\gamma-1}
		\bigl((P_R^0)_{jj}+(P_R^1)_{jj}\bigr)\\
		&=
		\int_{\mathbb R^d}
		(V_R(x)-V(x))f(x)\,dx+
		\int_{\mathbb R^d}
		(V_R(x)-V(x-Re_1))f_R(x)\,dx\\
		&=
		\int_{\mathbb R^d}
		\bigl((f(x)+f_R(x))^{\frac{p}{p-1}}-f^{\frac{p}{p-1}}(x)-f_R^{\frac{p}{p-1}}(x)\bigr)\,dx=A_R.
	\end{aligned}
\end{equation}
Here we used \eqref{eq:translated-functions},
\eqref{eq:trial-potential}, and \eqref{eq:interaction-term},
together with the identities $V=f^{1/(p-1)}$
and $V(\cdot-Re_1)=f_R^{1/(p-1)}$.
Finally, since $P_R$ is positive semidefinite and
$|\lambda_j|^{\gamma-1}>0$ for $1\leq j\leq K$,
we infer from \eqref{eq:trace-identity} that
\begin{equation}\label{eq:positive-norm}
	\begin{aligned}
		\|P_R\|_{\mathrm{op}}
		\leq \operatorname{Tr}P_R
		=\sum_{j=1}^K
		\bigl((P_R^0)_{jj}+(P_R^1)_{jj}\bigr)
		&\leq
		\frac{1}{\min\limits_{1\leq j\leq K}|\lambda_j|^{\gamma-1}}
		\sum_{j=1}^K|\lambda_j|^{\gamma-1}
		\bigl((P_R^0)_{jj}+(P_R^1)_{jj}\bigr)\\
		&=
		\frac{pL^{(N)}_{\gamma,d,s}}
		{\gamma\min\limits_{1\leq j\leq K}|\lambda_j|^{\gamma-1}}
		A_R\leq CA_R,
	\end{aligned}
\end{equation}
where $C>0$ is independent of $R$.

Set
\begin{equation}\label{eq:reference-matrix}
	\Lambda:=\operatorname{diag}(\lambda_1,\ldots,\lambda_K),
\qquad
H_0:=
\begin{pmatrix}
	\Lambda & 0\\
	0 & \Lambda
\end{pmatrix}.
\end{equation}
Let $H_R$ be the $2K\times2K$ matrix associated with
$H_{V_R}$ and the ordered family
$(u_1,\ldots,u_K,u_{1,R},\ldots,u_{K,R})$, defined by
\[
H_R:=
\begin{pmatrix}
	(\langle u_j,H_{V_R}u_k\rangle)_{j,k=1}^K
	&
	(\langle u_j,H_{V_R}u_{k,R}\rangle)_{j,k=1}^K\\
	(\langle u_{j,R},H_{V_R}u_k\rangle)_{j,k=1}^K
	&
	(\langle u_{j,R},H_{V_R}u_{k,R}\rangle)_{j,k=1}^K
\end{pmatrix}.
\]
Here all inner products are taken in $L^2(\mathbb R^d)$.
These matrix elements are well-defined, since
Theorem~\ref{thm:structure} gives
$u_j\in H^{2s}(\mathbb R^d)$, their translates
have the same regularity, and
$V_R\in L^\infty(\mathbb R^d)$.
We first compute the diagonal blocks of $H_R$.
Since $H_Vu_k=\lambda_k u_k$ and
$H_{V_R}=H_V-(V_R-V)$, we have
\[
\begin{aligned}
	\langle u_j,H_{V_R}u_k\rangle
	&=\langle u_j,H_Vu_k\rangle
	-\int_{\mathbb R^d}(V_R(x)-V(x))u_j(x)u_k(x)\,dx\\
	&=\lambda_k\langle u_j,u_k\rangle-(P_R^0)_{jk}=\lambda_k\delta_{jk}-(P_R^0)_{jk}.
\end{aligned}
\]
Similarly, using the eigenvalue equation for $u_{k,R}$,
we also have
\[
\begin{aligned}
	\langle u_{j,R},H_{V_R}u_{k,R}\rangle
	&=\langle u_{j,R},H_{V(\cdot-Re_1)}u_{k,R}\rangle-\int_{\mathbb R^d}
	(V_R(x)-V(x-Re_1))u_{j,R}(x)u_{k,R}(x)\,dx\\
	&=\lambda_k\langle u_{j,R},u_{k,R}\rangle-(P_R^1)_{jk}=\lambda_k\delta_{jk}-(P_R^1)_{jk}.
\end{aligned}
\]
For the upper-right off-diagonal block, by \eqref{eq:trial-uniform-bound} and \eqref{eq:state-overlap}, there holds
\[
\begin{aligned}
	|\langle u_j,H_{V_R}u_{k,R}\rangle|
	&=\left|\lambda_k\langle u_j,u_{k,R}\rangle
	-\int_{\mathbb R^d}(V_R(x)-V(x-Re_1))
	u_j(x)u_{k,R}(x)\,dx\right|\\
	&\leq
	\bigl(|\lambda_k|+\|V_R\|_\infty+\|V\|_\infty\bigr)
	\int_{\mathbb R^d}|u_j(x)u_{k,R}(x)|\,dx\\
	&\leq CR^{-(d+2s)}.
\end{aligned}
\]
The lower-left off-diagonal block is the transpose
of the upper-right block, since
$H_{V_R}$ is self-adjoint and the eigenfunctions are real.
We then introduce
\[
B_R:=
\begin{pmatrix}
	0 & (\langle u_j,H_{V_R}u_{k,R}\rangle)_{j,k=1}^K\\
	(\langle u_{j,R},H_{V_R}u_k\rangle)_{j,k=1}^K & 0
\end{pmatrix}.
\]
Since $K$ is fixed, the preceding estimates imply that
\begin{equation}\label{eq:energy-matrix}
	H_R=H_0-P_R+B_R,\qquad
	\|B_R\|_{\mathrm{op}}=O(R^{-(d+2s)}).
\end{equation}

For sufficiently large $R$, the positive definiteness
of $G_R$ allows us to orthonormalize the family
$(u_1,\ldots,u_K,u_{1,R},\ldots,u_{K,R})$
using $G_R^{-1/2}$.
Indeed, the Gram matrix of the resulting family is
$G_R^{-1/2}G_RG_R^{-1/2}=I_{2K}$.
The corresponding energy matrix is therefore
\[
M_R:=G_R^{-1/2}H_RG_R^{-1/2}.
\]
By using \eqref{eq:Gram}, the spectral theorem, and
the Taylor expansion of $(1+t)^{-1/2}$ at $t=0$,
we obtain
\begin{equation}\label{eq:Gram-inverse-expansion}
	G_R^{-1/2}
=I_{2K}-\frac12E_R+O(R^{-2(d+2s)}).
\end{equation}
Moreover, from \eqref{eq:overlap-convergence}, \eqref{eq:reference-matrix},
\eqref{eq:positive-norm}, and
\eqref{eq:energy-matrix}, it follows that
\begin{equation}\label{eq:energy-uniform-bound}
	\|H_R\|_{\mathrm{op}}
\leq\|H_0\|_{\mathrm{op}}
+\|P_R\|_{\mathrm{op}}
+\|B_R\|_{\mathrm{op}}
=O(1)
\qquad \text{as }R\to\infty.
\end{equation}
Substituting \eqref{eq:Gram-inverse-expansion} into
$M_R=G_R^{-1/2}H_RG_R^{-1/2}$, we infer from \eqref{eq:Gram-error-bound} and \eqref{eq:energy-uniform-bound} that
\begin{equation}\label{eq:orthogonalized-energy-expansion}
	M_R=H_R-\frac12(E_RH_R+H_RE_R)
+O(R^{-2(d+2s)}).	
\end{equation}
Here and below, matrix-valued remainder estimates
are understood in the operator norm.
Combining \eqref{eq:orthogonalized-energy-expansion} with \eqref{eq:energy-matrix},
we find
\[
\begin{aligned}
	M_R-H_0
	&=-P_R+B_R
	-\frac12\left(E_R(H_0-P_R+B_R)
	+(H_0-P_R+B_R)E_R\right)+O(R^{-2(d+2s)})\\
	&=-P_R+B_R-\frac12(E_RH_0+H_0E_R)+\frac12(E_RP_R+P_RE_R)
	-\frac12(E_RB_R+B_RE_R)\\
	&\quad +O(R^{-2(d+2s)}).
\end{aligned}
\]
Observe that \eqref{eq:Gram-error-bound} and
\eqref{eq:positive-norm} give
\[
\|E_RP_R\|_{\mathrm{op}}
+\|P_RE_R\|_{\mathrm{op}}
\leq2\|E_R\|_{\mathrm{op}}\|P_R\|_{\mathrm{op}}
=O(A_RR^{-(d+2s)}).
\]
Similarly, \eqref{eq:Gram-error-bound} and \eqref{eq:energy-matrix} yield
\[
\|E_RB_R\|_{\mathrm{op}}
+\|B_RE_R\|_{\mathrm{op}}
\leq2\|E_R\|_{\mathrm{op}}\|B_R\|_{\mathrm{op}}
=O(R^{-2(d+2s)}).
\]
Consequently, we arrive at
\begin{equation}\label{eq:matrix-expansion}
M_R-H_0
=-P_R+B_R-\frac12(E_RH_0+H_0E_R)+O\bigl(A_RR^{-(d+2s)}
+R^{-2(d+2s)}\bigr).		
\end{equation}
Moreover, since $H_0$ defined in \eqref{eq:reference-matrix}
is fixed, we derive from \eqref{eq:Gram-error-bound} that
\[
\|E_RH_0+H_0E_R\|_{\mathrm{op}}
\leq2\|H_0\|_{\mathrm{op}}\|E_R\|_{\mathrm{op}}
=O(R^{-(d+2s)}).
\]
Taking the operator norm in \eqref{eq:matrix-expansion},
we therefore obtain
\begin{equation}\label{eq:matrix-norm-estimate}
	\begin{aligned}
		\|M_R-H_0\|_{\mathrm{op}}
		&\leq\|P_R\|_{\mathrm{op}}
		+\|B_R\|_{\mathrm{op}}
		+\frac12\|E_RH_0+H_0E_R\|_{\mathrm{op}}+C\bigl(A_RR^{-(d+2s)}+R^{-2(d+2s)}\bigr)\\
		&=O\bigl(A_R+R^{-(d+2s)}\bigr).
	\end{aligned}
\end{equation}

We now derive the desired trace estimate.
Set $F(t)=(-t)^\gamma$ for $t<0$.
Since the spectrum of $H_0$ is bounded away from zero,
the standard finite-dimensional trace expansion gives
\begin{equation}\label{eq:trace-Taylor}
	\Tr F(H_0+Z)
	=\Tr F(H_0)+\Tr(F'(H_0)Z)
	+O(\|Z\|_{\mathrm{op}}^2)
\end{equation}
for every self-adjoint matrix $Z$ of sufficiently small
operator norm.
We apply \eqref{eq:trace-Taylor} with $Z=M_R-H_0$.
The matrix $F'(H_0)$ is block diagonal, whereas $B_R$ and
$E_RH_0+H_0E_R$ have zero diagonal blocks.
Consequently,
\begin{equation}\label{eq:off-diagonal-trace-vanishing}
	\Tr(F'(H_0)B_R)=0,
\qquad
\Tr\bigl(F'(H_0)(E_RH_0+H_0E_R)\bigr)=0.
\end{equation}
In addition, noting that $F'(\lambda_j)=-\gamma|\lambda_j|^{\gamma-1}$,
we deduce from \eqref{eq:trace-identity} that
\begin{equation}\label{eq:linear-gain}
	\begin{aligned}
		-\Tr(F'(H_0)P_R)
		&=-\Tr(F'(\Lambda)P_R^0)-\Tr(F'(\Lambda)P_R^1)\\
		&=\gamma\sum_{j=1}^K|\lambda_j|^{\gamma-1}
		\bigl((P_R^0)_{jj}+(P_R^1)_{jj}\bigr)=pL^{(N)}_{\gamma,d,s}A_R.
	\end{aligned}
\end{equation}
Furthermore, since $V$ is a normalized optimizer, we also have
\begin{equation}\label{eq:reference-trace}
	\Tr F(H_0)
=
2\sum_{j=1}^K|\lambda_j|^\gamma
=
2L^{(N)}_{\gamma,d,s}.
\end{equation}
Combining \eqref{eq:matrix-expansion}, \eqref{eq:matrix-norm-estimate}, 
\eqref{eq:trace-Taylor}, \eqref{eq:off-diagonal-trace-vanishing}, \eqref{eq:linear-gain}, and \eqref{eq:reference-trace}, we conclude that
\begin{equation}\label{eq:numerator}
	\begin{aligned}
		\Tr F(M_R)
		&=
		2L^{(N)}_{\gamma,d,s}
		+pL^{(N)}_{\gamma,d,s}A_R
		+O\bigl(A_R^2+R^{-2(d+2s)}\bigr)\\
		&=
		2L^{(N)}_{\gamma,d,s}
		+pL^{(N)}_{\gamma,d,s}A_R
		+o(A_R).
	\end{aligned}
\end{equation}
The last equality follows from $A_R\to0$ as $R\to\infty$
and \eqref{eq:overlap-small-o}.

\textbf{Step 3. Strict binding inequality and compactness.}
By \eqref{eq:matrix-norm-estimate} and
\eqref{eq:overlap-convergence}, we have
$\|M_R-H_0\|_{\mathrm{op}}\to0$ as $R\to\infty$.
Since $H_0\leq\lambda_K I_{2K}$ and $\lambda_K<0$,
it follows that
\[
M_R\leq H_0+\|M_R-H_0\|_{\mathrm{op}}I_{2K}
\leq\frac{\lambda_K}{2}I_{2K}<0
\]
for all sufficiently large $R$.
Denote the eigenvalues of $M_R$, arranged in
nondecreasing order and counted with multiplicity, by
\[
\lambda_1(M_R)\leq\cdots\leq\lambda_{2K}(M_R).
\]
Since $M_R$ is the energy matrix of $H_{V_R}$
in an orthonormal basis of the $2K$-dimensional space
spanned by $u_1,\ldots,u_K,u_{1,R},\ldots,u_{K,R}$,
the min--max principle gives
\[
\lambda_j(H_{V_R})\leq\lambda_j(M_R)<0,
\qquad 1\leq j\leq2K.
\]
Thus $H_{V_R}$ has at least $2K$ negative eigenvalues,
counted with multiplicity.
Since $\gamma>0$, we conclude that
\[
\sum_{j=1}^{2K}|\lambda_j(H_{V_R})|^\gamma
\geq\sum_{j=1}^{2K}|\lambda_j(M_R)|^\gamma
=\Tr F(M_R).
\]
Using $V_R$ as a trial potential in \eqref{eq:constant}
with $N$ replaced by $2K$, and recalling
\eqref{eq:denominator} and \eqref{eq:numerator},
we arrive at
\begin{equation}\label{eq:strict-quotient}
	\begin{aligned}
		L^{(2K)}_{\gamma,d,s}
		\geq
		\frac{\displaystyle\sum_{j=1}^{2K}
			|\lambda_j(H_{V_R})|^\gamma}{\|V_R\|_p^p}
		\geq\frac{\Tr F(M_R)}{2+A_R}
		&=\frac{2L^{(N)}_{\gamma,d,s}
			+pL^{(N)}_{\gamma,d,s}A_R+o(A_R)}{2+A_R}\\
		&=L^{(N)}_{\gamma,d,s}
		+\frac{(p-1)L^{(N)}_{\gamma,d,s}A_R+o(A_R)}
		{2+A_R}\\
		&>L^{(N)}_{\gamma,d,s}
	\end{aligned}
\end{equation}
for all sufficiently large $R$.
The last inequality follows from $p>1$,
$L^{(N)}_{\gamma,d,s}>0$, $A_R>0$, and
$o(A_R)/A_R\to0$ as $R\to\infty$.
Since $K\leq N$, the rank monotonicity
\eqref{eq:monotonicity} and \eqref{eq:strict-quotient}
imply that
\[
L^{(2N)}_{\gamma,d,s}
\geq L^{(2K)}_{\gamma,d,s}
>L^{(N)}_{\gamma,d,s}.
\]
As $N\in\mathbb N$ was arbitrary,
\eqref{eq:doubling} follows. In particular, 
\[
L_{\gamma,d,s}
\geq L^{(2N)}_{\gamma,d,s}
>L^{(N)}_{\gamma,d,s},
\qquad N\in\mathbb N.
\]

Finally, we establish the compactness assertion.
For $N=1$, the proof of Theorem~\ref{thm:existence}
already excludes vanishing, while dichotomy is impossible
since it would require an integer $1\leq k\leq N-1=0$.
Thus, after suitable translations, every normalized
maximizing sequence is tight, and
Corollary~\ref{cor:tight-maximizing} yields
precompactness modulo translations in
$L^p(\mathbb R^d)$.
For $N\geq2$, set $r=\lfloor N/2\rfloor$.
Since $1\leq r$ and $2r\leq N$, it follows from
\eqref{eq:doubling} and \eqref{eq:monotonicity} that
\[
\max_{1\leq k\leq\lfloor N/2\rfloor}
L^{(k)}_{\gamma,d,s}
=
L^{(r)}_{\gamma,d,s}
<
L^{(2r)}_{\gamma,d,s}
\leq
L^{(N)}_{\gamma,d,s}.
\]
Hence \eqref{eq:binding-criterion} holds, and Remark~\ref{rem:binding-compactness} indicates that every normalized maximizing sequence for $L^{(N)}_{\gamma,d,s}$ is precompact modulo translations in $L^p(\mathbb R^d)$.
This completes the proof of Theorem~\ref{thm:doubling}.
\end{proof}

{\bf Acknowledgments} 
L. Chen was supported by the National Natural Science Foundation of China (No.~12622108) and a grant from the Beijing Institute of Technology (No.~2022CX01002).
The authors acknowledge the use of AI tools for valuable assistance during the exploratory stage of this project. All mathematical arguments and proofs in the final manuscript were independently checked and written by the authors.

\smallskip

{\bf Data availability statement} Data sharing not applicable to this article as no datasets were generated or analysed during the current study.

\smallskip

{\bf Conflict of interest} The authors declare that they have no conflict of interest.

\end{document}